\documentclass[a4paper,oneside]{amsart}

\usepackage{amsmath}
\usepackage{amssymb}
\usepackage{amsthm}
\usepackage{tikz}
\usepackage{xypic}
\usepackage[pagebackref]{hyperref}

\newcommand{\Bigast}{\mathop{\scalebox{1.5}{\raisebox{-0.2ex}{$\ast$}}}}

\usetikzlibrary{decorations.markings}
\tikzstyle arrowstyle=[scale=1]
\tikzstyle directed=[postaction={decorate,decoration={markings,
    mark=at position .65 with {\arrow[arrowstyle]{stealth}}}}]

\DeclareMathOperator{\rank}{rank}
\DeclareMathOperator{\cone}{cone}

\DeclareMathOperator{\Tor}{Tor}
\DeclareMathOperator{\Ker}{Ker}
\DeclareMathOperator{\Img}{Im}
\DeclareMathOperator{\ab}{ab}
\DeclareMathOperator{\st}{st}
\DeclareMathOperator{\sk}{sk}
\DeclareMathOperator{\defi}{def}
\DeclareMathOperator{\supp}{supp}
\DeclareMathOperator{\conv}{conv}
\DeclareMathOperator{\Cart}{Cart}
\DeclareMathOperator{\Gen}{Gen}
\DeclareMathOperator{\Red}{Red}

\def\pt{{\mathrm{pt}}}

\def\cc{{\mathrm{cc}}}
\def\ZZ{\mathbb{Z}}
\def\QQ{\mathbb{Q}}

\def\RC{\mathrm{RC}}
\def\R{\mathcal{R}}
\def\Z{\mathcal{Z}}
\def\K{\mathcal{K}}
\def\L{\mathcal{L}}
\def\ZK{\Z_\K}
\def\RK{\R_\K}
\def\ccK{\cc(\K)}
\def\b{\widetilde{b}}
\def\G{\underline{G}}
\def\g{{\underline{g}}}
\def\h{{\underline{h}}}
\def\H{{\widetilde{H}}}
\def\k{\mathbf{k}}

\newtheorem*{thm*}{Theorem}
\newtheorem{thm}{Theorem}[section]
\newtheorem{lmm}[thm]{Lemma}

\newtheorem{prp}[thm]{Proposition}
\newtheorem{crl}[thm]{Corollary}

\theoremstyle{definition}

\newtheorem{dfn}[thm]{Definition}
\newtheorem{rmk}[thm]{Remark}
\newtheorem{exm}[thm]{Example}
\newtheorem{con}[thm]{Construction}
\newtheorem{prb}[thm]{Problem}
\numberwithin{equation}{section}

\hypersetup{
    colorlinks=true,
    linkcolor=blue,
    citecolor=blue,
    urlcolor=blue,
    pdftitle={Cartesian subgroups in graph products of groups},
    pdfauthor={Fedor Vylegzhanin},
}



\author{Fedor Vylegzhanin}
\address{\parbox{\linewidth}{
Steklov Mathematical Insitute of Russian Academy of Sciences, Moscow, Russia;\\
National Research University Higher School of Economics, Moscow, Russia.
}}
\email{vylegf@gmail.com}
\title[Cartesian subgroups in graph products]{Cartesian subgroups in graph products of groups}

\begin{document}
\begin{abstract}
The kernel of the natural projection of a graph product of groups onto their direct product is called the Cartesian subgroup of the graph product. This construction generalises commutator subgroups of right-angled Coxeter and Artin groups. Using theory of polyhedral products, we give a lower and an upper bound on the number of relations in presentations of Cartesian groups and on their deficiency. The bounds are related to the fundamental groups of full subcomplexes in the clique complex, and the lower bound coincide with the upper bound if these fundamental groups are free or free abelian. 

Following Li Cai's approach, we also describe an algorithm that computes ``small'' presentations of Cartesian subgroups.
\end{abstract}
\maketitle

\section{Introduction}
Let $\Gamma$ be a simple graph on the vertex set $[m]=\{1,\dots,m\}$ and $\G=(G_1,\dots,G_m)$ be a sequence of discrete groups. The corresponding \emph{graph product} \cite{green}
$$\G^\Gamma:=(G_1\ast\dots\ast G_m)/(g_ig_j=g_jg_i,~\forall\{i,j\}\in\Gamma,~\forall g_i\in G_i,~\forall g_j\in G_j)$$
interpolates between the free product and the direct product of $G_1,\dots,G_m$ as $\Gamma$ varies between the edgeless graph and the complete graph. The \emph{Cartesian subgroup}
$$\Cart(\G,\Gamma):=\Ker(\G^\Gamma\to G_1\times\dots\times G_m)$$
of the graph product was
studied in \cite{holt_rees,pv,pv_artin}.
As a group, $\Cart(\G,\Gamma)$ depends only on the graph $\Gamma$ and on cardinalities $|G_i|$ of the groups $G_i.$ This is well known for the classical Cartesian subgroup $\Ker(\Bigast_i G_i\to\prod_i G_i)$ which corresponds to the edgeless graph and is free on $\sum_{J\subset[m]}(|J|-1)\prod_{i\in J} |G_i\setminus\{1_i\}|$ generators \cite[Theorem 5.1]{gruenberg}. 

The  \emph{right-angled Coxeter groups} $$\RC_\Gamma:=\langle g_1,\dots,g_m\mid g_i^2=1,~i=1,\dots,m;~g_ig_j=g_jg_i,~\{i,j\}\in\Gamma\rangle = (\underline{\ZZ_2})^\Gamma$$ and \emph{right-angled Artin groups}
$\mathrm{RA}_\Gamma:=\langle g_1,\dots,g_m\mid g_ig_j=g_jg_i,~\{i,j\}\in\Gamma\rangle = \underline{\ZZ}^\Gamma$ are special cases of the graph product construction. Their Cartesian subgroups are the commutator subgroups $\RC_\Gamma'=\Cart(\underline{\ZZ_2},\Gamma)$ and $\mathrm{RA}_\Gamma'=\Cart(\underline{\ZZ},\Gamma).$ The classifying spaces of graph products and their Cartesian subgroups can be described as certain \emph{polyhedral products} of topological spaces \cite{pv,bbc} (see Proposition \ref{prp:polyhedral_products_are_classifying_spaces}); conversely, the fundamental groups of some polyhedral products can be identified with $\G^\Gamma$ or $\Cart(\G,\Gamma)$. For example, each simplicial complex $\K$ corresponds to the \emph{real moment-angle complex} $\RK$, a polyhedral product important in toric topology \cite{ToricTopology}. It is known that $\pi_1(\RK)\cong\RC'_{\Gamma}=\Cart(\underline{\ZZ_2},\Gamma)$ where $\Gamma=\sk_1\K$ is the $1$-skeleton of $\K$. (In general, fundamental groups of polyhedral products are \emph{relative graph products}, see \cite[Theorem 2.18]{davis_asphericity}.)

We study presentations of Cartesian subgroups and related numerical invariants: their \emph{rank} and \emph{deficiency}
$$\rank G := \inf_{G=\langle X\mid R\rangle}|X|,\quad \defi G := \inf_{G=\langle X\mid R\rangle} (|R|-|X|).$$

It is convenient to consider the \emph{clique complex} of the graph $\Gamma$
$$\K(\Gamma):=\{I\subset [m]:~\{i,j\}\in\Gamma,~\forall i,j\in I\},$$  the unique flag simplicial complex $\K$ such that $\Gamma=\sk_1\K.$ \textbf{From now on we assume that $\K=\K(\Gamma)$ and write $\G^\K:=\G^{\Gamma},$ etc, by abuse of notation.}

The ranks of Cartesian subgroups were computed by Panov and Veryovkin \cite[Theorem 5.2(b)]{pv_artin}: $\rank\Cart(\G,\K)=N(\G,\K),$ where
$$N(\G,\K):=\sum_{J\subset[m]}n_J\cdot\b_0(\K_J),$$
$n_J:=\prod_{j\in J}|G_j\setminus\{1_j\}|$ for $J\subset[m]$ and $\b_i(X):=\dim_\QQ\H_i(X;\QQ).$ (In particular, the induced subgraph $\Gamma_J$ of $\Gamma$ has $\b_0(\K_J)+1$ path components.) Moreover, Panov and Veryovkin described an explicit minimal generating set for $\Cart(\G,\K)$ that consists of nested iterated commutators \cite[Theorem 5.2(a)]{pv_artin}. A special case of this result, the formula
$\rank\RC_\K'=\sum_{J\subset[m]}\b_0(\K_J)$ and a set of $\sum_{J\subset[m]}\b_0(\K_J)$ generators for the group $\RC_\K'$, was obtained earlier by the same authors \cite[Theorem 4.5]{pv}.

We provide an alternative minimal set of generators for $\Cart(\G,\K),$ following the approach of Li Cai \cite{licai_slides}, and describe a small sufficient set of relations between them. Let $X=\bigsqcup_{\alpha}X_\alpha$ be a topological space with path components $\{X_\alpha\}$. Consider the group $$\Pi_1(X):=\Bigast_\alpha\pi_1(X_\alpha).$$
\begin{thm}
\label{thm:presentation_exists_intro}
Let $\K$ be a flag simplicial complex on the vertex set $[m]$ and $\G=(G_1,\dots,G_m)$ be a sequence of discrete groups. Then the group $\Cart(\G,\K)=\Ker(\G^\K\to\prod_{i=1}^m G_i)$ admits a presentation by $N(\G,\K)=\sum_{J\subset[m]}n_J\cdot\b_0(\K_J)$ generators modulo
$$M^+(\G,\K):=\sum_{J\subset[m]}n_J\cdot\rank\Pi_1(\K_J)$$
relations, where $n_J:=\prod_{j\in J}|G_j\setminus\{1_j\}|$.
In particular,
$$\defi\Cart(\G,\K)\leq M^+(\G,\K)-N(\G,\K).$$
\end{thm}
See Theorem \ref{thm:explicit_presentation} for a detailed description of this presentation.
The generators are explicit (each generator is of the form $L_\g(i,J):=\prod_{j\in J}g_j\cdot g_i^{-1}\cdot (\prod_{j\in J\setminus i}g_j)^{-1}$), while the relations can be computed by a recursive algorithm. Moreover, the generators and relations have a geometric meaning in terms of paths in certain polyhedral products, and this result can be extended to the case of infinite graph products (see Subsection \ref{subsec:infinite}).
We also give a homological lower bound on the size of a presentation:
\begin{thm}
\label{thm:lower_bound}
Any presentation of the group $\Cart(\G,\K)$ contains at least
$$M^-(\G,\K):=\rank\bigoplus_{J\subset[m]}\H_1(\K_J;\ZZ)^{\oplus n_J}$$
relations. Moreover, $\defi \Cart(\G,\K)\geq M^-(\G,\K)-N(\G,\K).$
\end{thm}
Since $(\Pi_1(X))_{\ab}=\H_1(X;\ZZ)$, the numbers $M^-(\G,\K)$ and $M^+(\G,\K)$ coincide if all full subcomplexes of $\K$ have free or free abelian fundamental groups. This is the case when $m$ is small enough (see Remark \ref{rmk:m_small}) or if $\K$ is a part of a surface triangulation (see Proposition \ref{prp:flag_surface_triangulation}). In these cases, we know both the rank and deficiency of Cartesian subgroups.

As a special case we obtain a small presentation of the commutator subgroup $\RC_\K'$ in any right-angled Coxeter group, and a lower bound on the number of relations in any of its presentations:
\begin{crl}
\label{crl:presentations_for_rck}
Let $\K$ be a flag simplicial complex on the vertex set $[m]$. Then
\begin{enumerate}
\item The commutator subgroup $\RC_\K'$ of the corresponding right-angled Coxeter group $\RC_\K$ admits a presentation by $\sum_{J\subset[m]}\b_0(\K_J)$ generators modulo $\sum_{J\subset[m]}\rank\Pi_1(\K_J)$ relations.

\item Any presentation of the group $\RC_\K'$ contains at least\\$\rank\bigoplus_{J\subset[m]}\H_1(\K_J;\ZZ)$ relations.

\item The deficiency $\defi\RC_\K'$ satisfies the inequalities
$$
\rank\bigoplus_{J\subset[m]}\H_1(\K_J;\ZZ)
\leq
\defi\RC_\K'+\sum_{J\subset[m]}\b_0(\K_J)
\leq
\sum_{J\subset[m]}\rank\Pi_1(\K_J).
\qed$$
\end{enumerate}

\end{crl}

Two earlier results can be deduced from the theorems above:
\begin{itemize}
\item Characterisation of simplicial complexes $\K$ such that $\Cart(\G,\K)$ is a free group \cite[Theorem 3.2]{holt_rees} (independently obtained in \cite[Theorem 4.3]{pv});
\item Characterisation of simplicial complexes $\K$ such that $\RC_\K'$ is a one-relator group \cite[Theorem 3.2(a)]{onerelator}.
\end{itemize}

\subsection*{Structure of paper}
In Section \ref{section:preliminaries} we discuss necessary definitions and known results about simplicial complexes, graph products of groups and polyhedral products of topological spaces. In Section \ref{section:lower_bound} we prove the lower bound (Theorem \ref{thm:lower_bound}). In Section \ref{section:presentation_description} we describe an explicit set of generators for $\Cart(\G,\K)$ and an algorithm that computes a presentation on these generators (see Theorem \ref{thm:explicit_presentation} which implies Theorem \ref{thm:presentation_exists_intro}).  In Section \ref{section:upper_bound} we prove Theorem \ref{thm:explicit_presentation}, applying a version of the van Kampen theorem to a certain explicit classifying space of $\Cart(\G,\K)$. Finally, in section \ref{section:discussion} we discuss possible generalisations and discuss the similarities between results of this paper and known results on graph products of algebras.

\subsection*{Acknowledgements}
The author woud like to thank Taras E. Panov for suggesting the problem and valuable advice, Li Cai, Anton A. Klyachko, Temurbek Rahmatullaev and Andrei Yu. Vesnin for helpful discussions, and the anonymous referee for careful reading of the text and important corrections.

\section{Preliminaries}
\label{section:preliminaries}
\subsection{Simplicial complexes}
A \emph{simplicial complex} $\K$ on the vertex set $V$ is a non-empty collection of finite subsets $I\subset V,$ called the \emph{faces}, that is closed under inclusion (i.e. if $J\subset I\in\K$, then $J\in\K$). Usually $V=[m]:=\{1,\dots,m\}.$ If $i\in V$ and $\{i\}\notin\K$, then $i$ is called a \emph{ghost vertex}.
We consider only
complexes without ghost vertices,
so we assume that
$\{i\}\in\K$ for every $i\in V.$

Given $J\subset V,$ the simplicial complex $\K_J:=\{I\in\K:~I\subset J\}$ on vertex set $J$ is called a \emph{full subcomplex} of $\K.$
\newpage

A \emph{missing face} of $\K$ is a set $J\subset V$ such that $J\notin\K$, but every proper subset of $J$ is a face of $\K.$ A complex $\K$ is \emph{flag} if all of its missing faces consist of two vertices; equivalently, a flag complex is the clique complex of its 1-skeleton ($I\in\K$ if and only if $\{i,j\}\in\sk_1\K$ for all $i,j\in I$). Every full subcomplex of a flag complex is flag.

Every face $I\in\K$ corresponds to a \emph{geometric face} $|I|=\mathrm{conv}\{e_i:~i\in I\}$ of the standard geometric simplex $\Delta^{m-1}=\mathrm{conv}\{e_i:i\in[m]\}\subset\mathbb{R}^m.$ The union of these geometric faces is called the \emph{geometric realization} $|\K|$ of a complex $\K.$

\subsection{Polyhedral products}
Let $(\underline{X},\underline{A})=\{(X_i,A_i)\}_{i=1}^m$ be a sequence of pairs of topological spaces and $\K$ be a simplicial complex on $[m].$ The corresponding \emph{polyhedral product} is the following subspace of $X_1\times\dots\times X_m$:
$$(\underline{X},\underline{A})^\K:=\bigcup_{I\in \K}\Big(\prod_{i\in I}X_i\times \prod_{i\in[m]\setminus I}A_i\Big).$$
It is clear that $(\underline{X},\underline{A})^\K\simeq (\underline{Y},\underline{B})^\K$ if $(X_i,A_i)\simeq (Y_i,B_i)$ for all $i=1,\dots,m$ (the symbol $\simeq$ denotes homotopy equivalence). We write $(X,A)^\K:=(\underline{X},\underline{A})^\K$ if $X_1=\dots=X_m=X,$ $A_1=\dots=A_m=A.$ For a sequence of pointed spaces $(X_1,\dots,X_m)$ denote also $\underline{X}^\K:=(\underline{X},\underline{\pt})^\K.$

The \emph{real moment-angle complex} $\RK:=(D^1,S^0)^\K$ is a special case of this construction. Homology groups of real moment-angle complexes are well known.
\begin{prp}[see {\cite[Theorem 4.5.8]{ToricTopology}}]\label{prp:hochster}
For any $i\geq 0$ and any group of coefficients, we have
$
    H_i(\RK)\cong\bigoplus_{J\subset[m]} \H_{i-1}(\K_J).
$
\qed
\end{prp}
(For $i=0$ we obtain $H_0(\RK;\k)\cong\H_{-1}(\varnothing;\k)\cong\k$.)
A far-reaching generalization of this result is the following stable homotopy decomposition for polyhedral products of the form
$(\cone \underline{A},\underline{A})^\K:$

\begin{thm}[{\cite[Theorem 2.21]{bbcg}}]
\label{thm:bbcg}
Let $(\underline{X},\underline{A})=\{(X_i,A_i)\}_{i=1}^m$ be a sequence of CW pairs such that all $X_i$ are contractible. Then there is a homotopy equivalence
$$
    \Sigma(\underline{X},\underline{A})^\K\simeq
    \Sigma^2\bigvee_{J\subset [m]}|\K_J|\wedge \underline{A}^{\wedge J},
$$
where $\underline{A}^{\wedge J}:=\bigwedge_{j\in J} A_j$ is the smash product.\qed
\end{thm}
\begin{crl}
\label{crl:bbcg}
Under the conditions of Theorem \ref{thm:bbcg}, we have, for any $i\geq 0$,
\[
	\H_i((\underline{X},\underline{A})^\K;\ZZ)\cong
	\bigoplus_{J\subset[m]}\H_{i-1}(|\K_J|\wedge\underline{A}^{\wedge J};\ZZ).
\qed
\]
\end{crl}

\subsection{Classifying spaces of graph products and Cartesian subgroups}
For a topological group $G$, there is a principal $G$-fibration $EG\to BG$ with $EG$ contractible, which is unique up to a weak homotopy equivalence. The base space $BG$ is called the \emph{classifying space} of $G.$ For discrete groups the classifying space $BG=K(G,1)$ is aspherical and $EG$ is its universal cover. In the next proposition we denote
$E\G:=(EG_1,\dots,EG_m)$ and $B\G:=(BG_1,\dots,BG_m).$

\begin{prp}[{see e.g. \cite[Proposition 3.1]{pv}}]
Let $\G=(G_1,\dots,G_m)$ be a sequence of topological groups and $\K$ be a simplicial complex on $[m]$. Then there is a canonical homotopy fibration
$(E\G,\G)^\K\to (B\G)^\K\to\prod_{i=1}^mBG_i.
\qed$
\end{prp}

Combining this fibration with the results of Panov, Ray and Vogt \cite[Proposition 5.1]{prv} provides an important connection between polyhedral products of spaces and graph products of discrete groups:
\begin{prp}[{\cite[Theorem 3.2]{pv}}]
\label{prp:polyhedral_products_are_classifying_spaces}
Let $\G=(G_1,\dots,G_m)$ be a sequence of discrete groups and $\K$ be a simplicial complex on $[m]$. Denote $\Gamma=\sk_1\K$. Then
\begin{enumerate}
\item $\pi_1((B\G)^\K)\cong \G^\Gamma$ and $\pi_1((E\G,\G)^\K)\cong\Cart(\G,\Gamma)$;
\item $\pi_k((B\G)^\K)\cong\pi_k((E\G,\G)^\K),~k\geq 2;$ 
\item $(E\G,\G)^\K$ and $(B\G)^\K$ are aspherical if and only if $\K$ is a flag complex.
\end{enumerate}
In particular, for flag complexes $\K$ we have
\[
B(\G^\K)=(B\G)^\K,\quad B(\Cart(\G,\K))=(E\G,\G)^\K.
\qed
\]
\end{prp}
\begin{rmk}
Davis described the fundamental groups $\pi_1((\underline{X},\underline{A})^\K)$ of general polyhedral products algebraically as \emph{relative graph products} \cite[Theorem 2.18]{davis_asphericity}. An almost complete answer to the question ``Which polyhedral products are aspherical?'' was obtained by Davis and Kropholler \cite[Theorem 1]{davis_asphericity_corrigenda}.

Asphericity of $(E\G,\G)^\K$ in the flag case can also be proved by applying the theory of CAT(0)-spaces to the universal covering over the cubical complex $(\cone\G,\G)^\K$, see \cite[Chapter 12]{davis} for $G_i=\ZZ_2$ case and \cite[Lemma 2.11]{davis_asphericity} in the general case. This approach can be traced back to Meier \cite[Section 4]{meier}.
\end{rmk}

Recall that the \emph{right-angled Coxeter group} $$\RC_\K:=\langle g_1,\dots, g_m\mid g_i^2=1,~i=1,\dots,m;~g_ig_j=g_jg_i,~\{i,j\}\in\K\rangle$$ is the graph product of groups $G_1=\dots=G_m=\ZZ_2.$

\begin{crl}[{\cite[Corollary 3.4]{pv}}]
\label{crl:RCK_RCK'_classifying}
Let $\K$ be a flag simplicial complex. Then
$$B(\RC_\K)=(\mathbb{R} P^\infty)^\K,\quad B(\RC_\K')=\RK.$$
\end{crl}
\begin{proof}
Indeed, by Proposition \ref{prp:polyhedral_products_are_classifying_spaces} we have $B(\RC_\K)=(B\ZZ_2)^\K=(\mathbb{R} P^\infty)^\K$ and $B(\RC'_\K)=(E\ZZ_2,\ZZ_2)^\K=(S^\infty,S^0)^\K\simeq (D^1,S^0)^\K=\RK$.
\end{proof}

The isomorphism $\pi_1(\RK)\cong\RC_\K'$ has a geometric interpretation: each generator $g_i\in\RC_\K$ corresponds to the path along the $i$th coordinate edge of the cube $(D^1)^m\supset \RK$. Group multiplication corresponds to concatenation of paths. Then the word $w\in\RC_\K$ belongs to the commutator subgroup if and only if the corresponding path is a loop. Every edge $\{i,j\}\in\K$ corresponds, from the algebraic point of view, to a relation $g_ig_j=g_jg_i$ in $\RC_\K,$ and, from the topological point of view, to $2^{m-2}$ squares in $\RK.$ Due to the presence of these squares, the paths $\dots g_ig_j\dots$ and $\dots g_jg_i\dots$ are homotopic in $\RK.$ Thus a homomorphism $\RC_\K'\to\pi_1(\RK)$ is well defined; in fact, it is an isomorophism. Its generalisation $\Cart(\G,\K)\overset\cong\longrightarrow \pi_1((\cone\G,\G)^\K)$ is discussed in Subsection \ref{subsec:loop_identification} below.

\subsection{Ranks of groups and free products}
The \emph{rank} of a group $G$ is the smallest number $N$ such that $G$ can be generated by $N$ elements. (If $G$ is not finitely generated, $N$ is considered as a  cardinal number.)

\begin{prp} The rank has the following properties:
\begin{enumerate}
\item (Grushko's theorem) $\rank\Bigast_\alpha G_\alpha= \sum_\alpha \rank G_\alpha;$
\item If $f:G\to H$ is a surjective homomorphism, then $\rank G\geq\rank H$;
\item $\rank\ZZ^m=m.$
\end{enumerate}
\end{prp}
\begin{proof}
\begin{enumerate}
\item This is \cite[Theorem 1.8]{lyndon}.
\item If $G$ is generated by $g_1,\dots,g_N,$ then $H$ is generated by $f(g_1),\dots,f(g_N).$
\item If $\ZZ^m$ could be generated by $n<m$ elements, then the vector space $\QQ^m$ would be a linear span of $n<m$ elements.
\qedhere
\end{enumerate}
\end{proof}

\begin{dfn}
Let $\{X_\alpha\}$ be the set of all path components of a topological space $X.$ Denote $\Pi_1(X):=\Bigast_\alpha \pi_1(X_\alpha).$
\end{dfn}
\begin{prp}
\label{prp:Pi1_properties}
Let $X$ be a topological space and $\{X_\alpha\}$ be the set of its path components. Then
\begin{enumerate}
\item $\rank\Pi_1(X)=\sum_\alpha\rank \pi_1(X_\alpha)$;
\item $\Pi_1(X)_{\ab}\cong H_1(X;\ZZ)$;
\item $\Pi_1(X)\cong\pi_1(\bigvee_\alpha X_\alpha)$.
\end{enumerate}
\end{prp}
\begin{proof}
Statement (1) follows from the Grushko's theorem,
(2) from the Poincar\'e--Hurewicz theorem \cite[Theorem 2A.1]{hatcher} and (3) from the van Kampen theorem \cite[Theorem 1.20]{hatcher}.
\end{proof}
\section{Proof of Theorem \ref{thm:lower_bound}}
\label{section:lower_bound}
\subsection{A general lower bound on the number of generators and relations}
Recall that the integer homology of a discrete group is isomorphic to the integer homology of its classifying space \cite[Proposition II.4.1]{brown}:
$$\mathrm{H}_i(G;\ZZ):=\Tor^{\ZZ[G]}_i(\ZZ,\ZZ)\cong H_i(BG;\ZZ).$$

The following lemma is essentially
\cite[\S II.5, Exercise 5a]{brown}. Notably, it has a purely algebraic proof \cite[Lemma 1.2]{epstein}.

\begin{lmm}
\label{lmm:low_est}
Let $G=\langle x_1,\dots,x_N\mid r_1,\dots,r_M\rangle$ be a finite group presentation. Then
$$N\geq \rank \mathrm{H}_1(G;\ZZ),\quad M-N\geq \rank \mathrm{H}_2(G;\ZZ)-\dim_\QQ \mathrm{H}_1(G;\QQ).
$$

In particular, $M\geq\rank \mathrm{H}_2(G;\ZZ).$
\end{lmm}
\begin{proof}
We construct a CW classifying space for the group $G$ by attaching cells to its presentation complex. Then
$
  BG = \bigvee_{i=1}^N S_i^1\cup \bigcup_{j=1}^M e_j^2\cup\bigcup_\beta e_{\beta}^{>2},
$ 
where the
2-cells are attached by the maps corresponding to the words $r_1,\dots,r_M\in F(x_1,\dots,x_N)\cong\pi_1(\bigvee_{i=1}^N S_i^1).$
The cellular chain complex of $BG$ is of the form $$0\leftarrow \ZZ\overset{~0}{\longleftarrow} \ZZ^N\overset{~\partial_2}{\longleftarrow} \ZZ^M \overset{~\partial_3}{\longleftarrow} \dots$$
and has $H_*(BG;\ZZ)\cong\mathrm{H}_*(G;\ZZ)$ as its homology.

Denote $k=\rank\partial_2$. Then $\mathrm{H}_1(G;\ZZ)$ is a quotient of $\ZZ^N$, $\mathrm{H}_2(G;\ZZ)$ is a quotient of $\Ker\partial_2\simeq\ZZ^{M-k}$, and $\mathrm{H}_1(G;\QQ)\simeq\QQ^{N-k}$. Hence $\rank\mathrm{H}_1(G;\ZZ)\leq N$ and $\rank\mathrm{H}_2(G;\ZZ)-\dim_\QQ\mathrm{H}_1(G;\QQ)\leq (M-k)-(N-k)=M-N$.
\end{proof}
\begin{prp}
\label{prp:rel_in_comm}
Let $G=\langle x_1,\dots,x_N\mid r_1,\dots,r_M\rangle$ be a finitely presented group such that $\mathrm{H}_1(G;\ZZ)\simeq\ZZ^N.$ Then $r_1,\dots,r_M$ belong to the commutator subgroup of $F(x_1,\dots,x_N).$
\end{prp}
\begin{proof} In the notation of the previous lemma, we have $\ZZ^N/\Img\partial_2\simeq \ZZ^N,$ hence $\partial_2=0.$ On the other hand, $\partial_2:\ZZ^M\to \ZZ^N$ is the abelianization of the map $F(r_1,\dots,r_M)\to F(x_1,\dots,x_N).$ Therefore $r_j\in \Ker(F(x_1,\dots,x_N)\overset{\ab}{\longrightarrow} \ZZ^N)=F(x_1,\dots,x_N)'.$
\end{proof}
\begin{rmk}
\label{rmk:low_est_for_algebras}
For a finitely generated abelian group $A,$ let $\mathrm{gen}\,A$ and $\mathrm{rel}\,A$ be the smallest numbers such that $A\simeq\ZZ^{\mathrm{gen}\,A}/\ZZ^{\mathrm{rel}\,A}.$ Clearly, $\mathrm{gen}\,A=\rank A$ and $\mathrm{rel}\,A=\rank A-\dim_\QQ(A\otimes_\ZZ\QQ).$ 
Hence Lemma \ref{lmm:low_est} can be stated as
$$N\geq\mathrm{gen}\,\mathrm{H}_1(G;\ZZ),\quad M-N\geq\mathrm{gen}\,\mathrm{H}_2(G;\ZZ)-\mathrm{gen}\,\mathrm{H}_1(G;\ZZ)+\mathrm{rel}\,\mathrm{H}_1(G;\ZZ).$$
Since $\mathrm{H}_i(G;\ZZ)=\Tor^{\ZZ[G]}_i(\ZZ,\ZZ),$ we deduce: if $G$ has a presentation by $N$ generators and $M$ relations, then
$$N\geq\mathrm{gen}\,\Tor^{\ZZ[G]}_1(\ZZ,\ZZ),~M\geq\mathrm{gen}\,\Tor^{\ZZ[G]}_2(\ZZ,\ZZ)+\mathrm{rel}\,\Tor^{\ZZ[G]}_1(\ZZ,\ZZ).$$
There is a similar result on presentations of connected graded associative algebras with unit \cite[Theorem A.10]{hozk_flag}: if $\k$ is a principal ideal domain and $R=\bigoplus_{n\geq 0}R_n$ is a connected graded $\k$-algebra presented by $N$ generators modulo $M$ relations, then
$$N\geq\mathrm{gen}\,\Tor^R_1(\k,\k),~M\geq\mathrm{gen}\,\Tor^R_2(\k,\k)+\mathrm{rel}\,\Tor^R_1(\k,\k).$$
Moreover, there is a ``minimal'' presentation of $R$ such that both lower bounds are achieved. This is not true for groups: if $G$ is a knot group, then $\mathrm{H}_1(G;\ZZ)=\ZZ,$ but usually $\rank G>1$. 
\end{rmk}
\subsection{Homology of Cartesian subgroups}
\begin{prp}
\label{prp:homology_of_cartgk}
Let $\G=(G_1,\dots,G_m)$ be a sequence of discrete groups and $\K$ be a flag simplicial complex on vertex set $[m]$. Then
$$\mathrm{H}_i(\Cart(\G,\K);\ZZ)\cong\bigoplus_{J\subset[m]}\H_{i-1}(\K_J;\ZZ)^{\oplus n_J},$$
where $n_J:=\prod_{j\in J}|G_j\setminus\{1_j\}|.$ In particular,
$\mathrm{H}_i(\RC_\K';\ZZ)\cong\bigoplus_{J\subset[m]}\H_{i-1}(\K_J;\ZZ).$
\end{prp}
\begin{proof}
We have $\mathrm{H}_*(\Cart(\G,\K);\ZZ)=H_*((E\G,\G)^\K;\ZZ)$ by Proposition \ref{prp:polyhedral_products_are_classifying_spaces}, and
$$H_i((E\G,\G)^\K;\ZZ)\cong \bigoplus_{J\subset[m]}\H_{i-1}(|\K_J|\wedge \G^{\wedge J};\ZZ)$$
by Corollary \ref{crl:bbcg}. Note that $G_j$ are discrete spaces. If $A$ and $B$ are discrete and $|A|=n+1,$ $|B|=m+1$ then $|A\wedge B|$ is discrete and has $mn+1$ points; hence $\G^{\wedge J}$ is discrete and has $n_J+1$ points. Finally, if $A$ is discrete and $|A|=n+1,$ then $X\wedge A\cong X^{\vee n},$ hence $\H_*(X\wedge A;\ZZ)\cong \H_*(X;\ZZ)^{\oplus n}.$
\end{proof}

\begin{proof}[{Proof of Theorem \ref{thm:lower_bound}}]
We have $$\rank \mathrm{H}_1(\Cart(\G,\K);\ZZ)=N(\G,\K),\quad\rank \mathrm{H}_2(\Cart(\G,\K);\ZZ)=M^-(\G,\K)$$
by Proposition \ref{prp:homology_of_cartgk}. Then, by Lemma \ref{lmm:low_est}, the deficiency of $\Cart(\G,\K)$ is not less than $M^-(\G,\K)-N(\G,\K)$, and any presentation of this group has at least $M^-(\G,\K)$ relations.
\end{proof}

We also obtain a curious property of ``minimal presentations'' of $\Cart(\G,\K)$, which is not obvious for the presentation that will be obtained in Theorem \ref{thm:explicit_presentation}.
\begin{prp}
\label{prp:relations_are_commutators}
Suppose that the number $N(\G,\K)$ is finite, and let $\Cart(\G,\K)=\langle x_1,\dots,x_{N(\G,\K)}\mid r_1,\dots,r_M\rangle$ be a group presentation. Then the relations $r_1,\dots,r_M$ belong to the commutator subgroup of $F(x_1,\dots,x_{N(\G,\K)}).$
\end{prp}
\begin{proof}
This follows from Proposition \ref{prp:rel_in_comm}, since $\mathrm{H}_1(\Cart(\G,\K);\ZZ)\simeq \ZZ^{N(\G,\K)}.$
\end{proof}

\section{Explicit presentations of Cartesian subgroups in graph products}\label{section:presentation_description}
In Theorem \ref{thm:explicit_presentation} we will describe a small presentation of the group $\Cart(\G,\K)$, clarifying Theorem \ref{thm:presentation_exists_intro}. The proof will be given in Section \ref{section:upper_bound}.

Denote by $F(X)$ the free group generated by a set $X$. We first describe a set of elements $D=\{L_\g(i,J):\dots\}\subset\Cart(\G,\K)$ and its subset $\widehat{D}\subset D$ of \emph{distinguished} elements. Then we show that each element of $D$ is equal in $\Cart(\G,\K)$ to a word on the distinguished elements, and provide an algorithm that computes such words. Hence each word $r\in F(D)$ is equal in $\Cart(\G,\K)$ to a word $\Red(r)\in F(\widehat{D})$ which can be computed inductively. Our presentation is of the form $$\Cart(\G,\K)=\langle \widehat{D}\mid \Red(R_\g(\lambda,J)):\dots\rangle,$$ where the words $R_\g(\lambda,J)\in F(\widehat{D})$ are known explicitly.

\subsection{Generators} Fix a sequence $\G=(G_1,\dots,G_m)$ of discrete groups and a flag simplicial complex $\K$ on $[m]$. Denote $G_j^*:=G_j\setminus\{1_j\}$.

\begin{dfn}
\label{dfn:LiJ}
For $i\in J\subset[m]$ and $\g=(g_j:j\in J)\in\prod_{j\in J}G_j^*,$ define
$$L_\g(i,J):=\prod_{j\in J}g_j \cdot g_i^{-1}\cdot \Big(\prod_{j\in J\setminus i}g_j\Big)^{-1}\in \Cart(\G,\K),$$
where the products are in the ascending order. 
Denote $$D:=\Big\{
L_\g(i,J):J\subset[m],~i\in J,~\g\in\prod_{j\in J}G_j^*
\Big\}\subset\Cart(\G,\K).$$
For example, $L_\g(3,\{1,2,3,8\})=g_1g_2g_3g_8\cdot g_3^{-1}\cdot g_8^{-1}g_2^{-1}g_1^{-1}.$
\end{dfn}
We use the following notation from \cite[Definition 5.2]{hozk_flag} inspired by \cite{gptw}.
\begin{dfn}
For each $J\subset[m]$, choose a subset $\Theta(J)\subset J\setminus\{\max(J)\}$ such that $\Theta(J)\sqcup\{\max(J)\}$ contains exactly one vertex from each path component of $|\K_J|$.
\end{dfn}
We define the set of \emph{distinguished} elements
$$\widehat{D}:=
\Big\{
L_\g(i,J):J\subset[m],~i\in\Theta(J),~\g\in\prod_{j\in J}G_j^*
\Big\}\subset D.$$
Since $|\Theta(J)|=\b_0(\K_J)=b_0(\K_J)-1$, we have $|\widehat{D}|=\sum_{J\subset[m]}n_J\b_0(\K_J)$ while $|D|=\sum_{J\subset[m]}n_J|J|.$ When $i\in\Theta(J)$, we write $\widehat{L}_\g(i,J)$ instead of $L_\g(i,J)$ to emphasize that $L_\g(i,J)\in\widehat{D}.$
\begin{rmk}
The following ``canonical'' choice of $\Theta(J)\subset J$ is used in \cite{gptw} and \cite{pv,pv_artin}: define $\Theta(J)$ to be the set of all $i\in J$ such that
\begin{enumerate}
\item The vertices $i$ and $\max(J)$ are in different path components of $|\K_J|;$
\item The vertex $i$ has the smallest number in its path component of $|\K_J|.$
\end{enumerate}
Then the Panov--Veryovkin set of generators for the group $\Cart(\G,\K)$ is naturally indexed by the same set. Indeed, for $i\in J\subset[m]$ and $\g\in\prod_{j\in J}G_j^*$ define $$\Gamma_\g(i,J) :=(g_{k_1},(g_{k_2},\dots(g_{k_s},(g_{\max(J)},g_i))\dots))\in\Cart(\G,\K),$$ 
where $J\setminus\{i,\max(J)\}=\{k_1<\dots<k_s\}$ and $(g,h):=g^{-1}h^{-1}gh$ is the group commutator. 
Then $\{\Gamma_\g(i,J):i\in\Theta(J),J\subset[m]\}$ is precisely the set of generators from \cite[Theorem 5.2]{pv_artin}.

Panov--Veryovkin generators are important for the analogy between the group $\pi_1(\RK)$ and the Hopf algebra $H_*(\Omega\ZK;\k)$ (see Subsection \ref{subsection:comparison} below). However, our choice of generators is more suitable for computations.
\end{rmk}
\begin{rmk} For right-angled Coxeter groups $\G^\K=(\underline{\ZZ_2})^\K=\RC_\K$ we omit the subscript $\g$ since then $|G_j^*|=1$. Li Cai proved (\cite{licai_slides}, unpublished) that the set $\widehat{D}=\{\widehat{L}(i,J):J\subset[m],i\in\Theta(J)\}\subset\RC_\K'$ generates the group $\RC_\K'.$ His approach relies on a study of the \emph{Davis complex} \cite[Chapter 7]{davis}, a contractible cubical complex which is the universal covering of $\RK.$
\end{rmk}

\subsection{Reduction to the distinguished elements}
\begin{lmm}
\label{lmm:elementary_relation}
Let $J\subset[m]$, $\{i,j\}\in\K_J$ and $\g\in \prod_{j\in J}G_j^*.$ Then
\begin{equation}
\label{eqn:elementary_relation}
    L_\g(i,J)L_\g(j,J\setminus i)=
    L_\g(j,J)L_\g(i,J\setminus j)\in\Cart(\G,\K).
\end{equation}
\end{lmm}
\begin{proof}
Indeed, we have $
L_\g(i,J)\cdot L_\g(j,J\setminus i)=
\prod_{k\in J}g_k\cdot g_i^{-1}g_j^{-1}\cdot (\prod_{k\in J\setminus\{i,j\}}g_k)^{-1},$
and the right hand side is symmetrical on $i,j$ since $g_i$ and $g_j$ commute.
\end{proof}

Also, $\prod_{j\in J}g_j=\prod_{j\in J\setminus\{\max(J)\}}g_j\cdot g_{\max(J)}$ and hence 
\begin{equation}
\label{eqn:boundary_relation}
L_\g(\max(J),J)=1\in\Cart(\G,\K).
\end{equation}

\begin{prb}
Apply the Reidemeister-Schreier algorithm \cite[Section II.4]{lyndon} to the subgroup $\Cart(\G,\K)\subset\G^\K$. Is it true that one obtains the presentation of $\Cart(\G,\K)$ by the generators $\{L_\g(i,J):i\in J\subset[m],~\g\in\prod_{j\in J}G_j^*\}$ modulo the relations \eqref{eqn:elementary_relation} and \eqref{eqn:boundary_relation}? This larger presentation would be explicit and natural with respect to maps of both groups and simplicial complexes. Probably, Theorem \ref{thm:explicit_presentation} can be proved by applying Tietze transformations to this presentation. However, this approach seems tedious.
\end{prb}

\begin{prp}
\label{prp:reduction_of_generators}
Every element $L_\g(i,J)\in D$ is equal in $\Cart(\G,\K)$ to a word  $\Red(L_\g(i,J))\in F(\widehat{D})$ on the distinguished elements. This word can be chosen naturally with respect to inclusions $\K\hookrightarrow\L$ of full subcomplexes.
\end{prp}
\begin{proof}
We say that $|J|$ is the \emph{length} of $L_\g(i,J)$, and argue by induction on length. The base cases $|J|\leq 2$: $L_\g(i,\{i\})=1;$ $L_\g(i,\{i,j\})=1$ if $\{i,j\}\in\K;$ otherwise $L_\g(i,\{i,j\})=\widehat{L}_\g(i,\{i,j\})$ for $i<j$ and $L_\g(i,\{i,j\})=\widehat{L}_\g(j,\{i,j\})^{-1}$ for $i>j.$

At the inductive step, for an element $L_\g(i,J)\in D$ there are two cases:
\begin{itemize}
\item $i$ and $\max(J)$ are in the same path component of $\K_J.$ Choose a path $$(i=i_0,i_1,\dots,i_k=\max(J))$$ in $\K_J,$ i.e., a sequence of edges $\{i_t,i_{t+1}\}\in\K_J$. Applying the relation \eqref{eqn:elementary_relation}, we replace $L_\g(i_t,J)$ with a product of $L_\g(i_{t+1},J)$ and generators of smaller length. On the last step we replace $L_\g(\max(J),J)$ with $1$ (using \eqref{eqn:boundary_relation}). Therefore, $L_\g(i,J)$ is replaced with a word on elements of smaller length. By the inductive hypothesis, each of this elements can be expressed through the distinguished elements.

\item $i$ and $\max(J)$ are in different path components of $\K_J.$ Then there is exactly one vertex $i'\in\Theta(J)$ such that $i$ and $i'$ are in the same path component. Choose a path $(i=i_0,\dots,i_k=i')$ in $\K_J.$ Arguing as above, we obtain a word on $L_\g(i',J)=\widehat{L}_\g(i',J)$ and on elements of smaller length.
\end{itemize}
To obtain functoriality, choose each time the lexicographically minimal path from $i$ to $i_k.$ This choice is preserved by inclusions of simplicial complexes.
\end{proof}
 \begin{rmk}
The algorithm of Proposition \ref{prp:reduction_of_generators} does not depend on choice of $\g\in\prod_jG_j^*$ and on the group structure in $G_1,\dots,G_m$. In more detail, if
$$\Red(L_\g(i,J))=\widehat{L}_\g(i_1, J_1)^{\varepsilon_1}\cdot\dotso\cdot \widehat{L}_\g(i_N,J_N)^{\varepsilon_N}$$
in $\Cart(\G,\K)$, then
$$\Red(L_\h(i,J))=\widehat{L}_\h(i_1, J_1)^{\varepsilon_1}\cdot\dotso\cdot \widehat{L}_\h(i_N,J_N)^{\varepsilon_N}$$
in $\Cart(\underline{H},\K)$
for any $\underline{H}=(H_1,\dots,H_m)$ and $\h\in\prod_jH_j^*.$
\end{rmk}

\subsection{Statement of the main theorem}
\label{subsection:statement_main}

\begin{dfn} Let $W\in F(D)$. Express the letters of $W$ through the distinguished elements using the algorithm from Proposition \ref{prp:reduction_of_generators}. The resulting word will be denoted by $\Red(W)\in F(\widehat{D}).$
\end{dfn}

\begin{dfn} Let $J\subset[m]$, $\g\in\prod_{j\in J}G_j^*$. Let $\lambda=(i_1,\dots,i_{k+1}=i_1)$ be a \emph{cycle} in $\K_J,$ i.e. a sequence of edges $\{i_t,i_{t+1}\}\in\K_J$ for $t=1,\dots,k.$ Define the word
$$R_\g(\lambda,J):= \prod_{t=1}^kL_\g(i_{t+1}, J\setminus i_{t})\cdot L_\g(i_{t}, J\setminus i_{t+1})^{-1}\in F(D).$$
\end{dfn}
For example, if $J=\{4,5,7,8\}$ and $\lambda=(7,5,4,7),$ we have
$$R_\g(\lambda,J)=
L_\g(5,458)L_\g(7,478)^{-1}\cdot L_\g(4,478)L_\g(5,578)^{-1}\cdot
L_\g(7,578)L_\g(4,458)^{-1}.$$

Recall that $\Pi_1(X):=\Bigast_\alpha\pi_1(X_\alpha),$
where $X=\bigsqcup_\alpha X_\alpha$ is the decomposition of $X$ onto its path components. By CW-approximation, each group $\Pi_1(\K_J)$ is generated by cycles in $\K_J$. 
\begin{thm}
\label{thm:explicit_presentation}
Let $\G=(G_1,\dots,G_m)$ be a sequence of discrete groups and $\K$ be a flag simplicial complex on vertex set $[m].$ For each $J\subset[m],$ choose a set $\Gen(J)$ of cycles $\lambda=(i_1,\dots,i_k,i_{k+1}=i_1)$ in $\K_J$ such that their images generate the group $\Pi_1(\K_J)$.

Then the group $\Cart(\G,\K)$ is presented by $\sum_{J\subset[m]}\b_0(\K_J)\cdot n_J$ generators
$$\widehat{D}=\{\widehat{L}_\g(i,J): J\subset[m],~\g\in\prod_{j\in J}G_j^*,~i\in\Theta(J)\}$$ modulo $\sum_{J\subset[m]}|\Gen(J)|\cdot n_J$ relations
$$\{\Red(R_\g(\lambda,J))=1: J\subset[m],\g\in \prod_{j\in J}G_j^*,~\lambda\in\Gen(J)\}.$$
\end{thm}
The proof will be given in section \ref{section:upper_bound}.
\begin{proof}[Proof of Theorem \ref{thm:presentation_exists_intro}]
By CW-approximation and Proposition \ref{prp:Pi1_properties}(1), for each $J\subset[m]$ there exists a minimal generating set $\Gen(J)$ for the group $\Pi_1(\K_J)$ represented by cycles, i.e. $|\Gen(J)|=\rank\Pi_1(\K_J)$. The corresponding presentation from Theorem \ref{thm:explicit_presentation} meets the conditions of Theorem \ref{thm:presentation_exists_intro}.
\end{proof}

\begin{rmk}
For $\g\in\prod_{j=1}^mG_j$, denote $\supp\g:=\{j\in[m]:~g_j\neq 1_j\}$. Then we can rewrite $\widehat{D}=\{\widehat{L}_\g(i,\supp\g):~\g\in\prod_{j=1}^m G_j,~i\in\Theta(\supp\g)\}$ and the set of relations as $\{\Red(R_\g(\lambda,\supp\g))=1:~\g\in\prod_{i=1}^m G_j,~\lambda\in\Gen(\supp\g)\}$.
\end{rmk}
\subsection{Examples}
\begin{exm}[$4$-cycle]
Let $\K$ be a four-cycle. Since $\K_{\{1,3\}}$ and $\K_{\{2,4\}}$ are the only disconnected full subcomplexes, $\Theta(\{1,3\})=\{1\},$ $\Theta(\{2,4\})=\{2\}$, all other $\Theta(J)$ being zero. We denote the corresponding distinguished elements
$\widehat{L}_\g(1,\{1,3\})=g_1g_3\cdot g_1^{-1}\cdot g_3^{-1}$ 
and
$\widehat{L}_\g(2,\{2,4\})=g_2g_4\cdot g_2^{-1}\cdot g_4^{-1}$ 
by $A_{g_1,g_3}$ and $B_{g_2,g_4}.$ It follows that the group
$\Cart(\G,\K)=\Ker((G_1\ast G_3)\times (G_2\ast G_4)\to G_1\times G_2\times G_3\times G_4)$
is generated by the set
$$\widehat{D}=\widehat{A}\sqcup\widehat{B}=
\{A_{g_1,g_3}:g_1\in G_1^*,g_3\in G_3^*\}\sqcup
\{B_{g_2,g_4}:g_2\in G_2^*,g_4\in G_4^*\}$$
modulo the following relations: for every $\g\in \prod_{j=1}^4 G_j^*,$ there is a relation\\$\Red(R_\g(\lambda,\{1,2,3,4\}))=1,$ where $\lambda=(1,2,3,4,1).$ The left hand side is equal to
\begin{align*}
\Red\Big(
&
L_\g(2,234)
L_\g(1,134)^{-1}
L_\g(3,134)
L_\g(2,124)^{-1}\\
\cdot&
L_\g(4,124)
L_\g(3,123)^{-1}
L_\g(1,123)
L_\g(4,234)^{-1}\Big)=\\
&\widehat{L}_\g(2,24)\cdot
\widehat{L}_\g(1,13)^{-1}
\cdot 1\cdot
\widehat{L}_\g(2,24)^{-1}
\cdot 1\cdot
1^{-1}\cdot
\widehat{L}_\g(1,13)\cdot
1^{-1}
=(B_{g_2,g_4}^{-1},A_{g_1,g_3}).
\end{align*}

Equivalently, $A_{g_1,g_3}$ and $B_{g_2,g_4}$ commute, hence $\Cart(\G,\K)\cong F(\widehat{A})\times F(\widehat{B})$ is a direct product of two free groups. In particular, $\Cart(\underline{\ZZ_2},\K)=\RC'_\K\cong\ZZ^2$ is generated by $A_{g_1,g_3}$ and $B_{g_2,g_4}.$
\end{exm}
\begin{exm}[$m$-cycles]
\label{exm:groups_cycles}
Let $\K_m$ be an $m$-cycle, $m>4.$ By \cite[Theorem 4.1.7]{ToricTopology}, $\RK$ is a closed $(n+1)$-manifold whenever $\K$ is a triangulation of $n$-sphere. The $m$-cycle corresponds to the oriented surface $S_{g(m)}\cong\R_{\K_m}$ of genus $g(m)=1+(m-4)2^{m-3},$ hence $\RC'_{\K_m}$ is the surface group $\pi_1(S_{g(m)}).$

Consider the presentation given by Theorem \ref{thm:explicit_presentation}.
The only full subcomplex with non-trivial $\Pi_1$ is the whole $\K_m.$ If the vertices are ordered cyclically, we have $\lambda=(1,2,\dots,m,1)$ and hence the following set of $\prod_{i=1}^m |G_i^*|$ defining relations for $\Cart(\G,\K_m):$
$$\Big\{\Red\Big(
\prod_{t=1}^m
L_\g(t+1,[m]\setminus\{t\})
L_\g(t,[m]\setminus\{t+1\})^{-1}
\Big)=1,~\forall\g\in\prod_{i=1}^m G_i^*\Big\}.
$$
In particular, $\RC'_{\K_m}$ is defined by the single relation $R_m=1,$ where
$$R_m:=\Red\Big(
\prod_{t=1}^m
L(t+1,[m]\setminus\{t\})
L(t,[m]\setminus\{t+1\})^{-1}
\Big).$$

For $m=5$ we have the relation
\begin{align*}\Red
\Big(
&
L(2,2345)
L(1,1345)^{-1}
\cdot
L(3,1345)
L(2,1245)^{-1}
\\
\cdot&
L(4,1245)
L(3,1235)^{-1}
\cdot
L(5,1235)
L(4,1234)^{-1}
\cdot
L(1,1234)
L(5,2345)^{-1}
\Big)=
\\
&
\widehat{L}(2,235)
\widehat{L}(3,35)
\widehat{L}(2,25)^{-1}
\widehat{L}(2,245)
\widehat{L}(3,35)^{-1}
\widehat{L}(1,134)^{-1}\\
\cdot&
\widehat{L}(3,135)
\widehat{L}(1,14)
\widehat{L}(2,245)^{-1}
\widehat{L}(1,124)^{-1}
\widehat{L}(2,25)
\widehat{L}(3,135)^{-1}
\widehat{L}(1,13)\\
\cdot&
\widehat{L}(2,235)^{-1}
\widehat{L}(1,13)^{-1}
\widehat{L}(1,124)
\widehat{L}(2,24)
\widehat{L}(1,14)^{-1}
\widehat{L}(1,134)
\widehat{L}(2,24)^{-1}.
\end{align*}
\begin{figure}[h!]
\begin{tikzpicture}[scale=0.5]
  
 \draw[thick] circle(5);

\draw[fill] (  0:5) circle (.1);
\draw[fill] ( 18:5) circle (.1);
\draw[fill] ( 36:5) circle (.1);
\draw[fill] ( 54:5) circle (.1);
\draw[fill] ( 72:5) circle (.1);
\draw[fill] ( 90:5) circle (.1);
\draw[fill] (108:5) circle (.1);
\draw[fill] (126:5) circle (.1);
\draw[fill] (144:5) circle (.1);
\draw[fill] (162:5) circle (.15);
\draw[fill] (180:5) circle (.1);  
\draw[fill] (198:5) circle (.1);
\draw[fill] (216:5) circle (.1);
\draw[fill] (234:5) circle (.1);
\draw[fill] (252:5) circle (.1);
\draw[fill] (270:5) circle (.1);
\draw[fill] (288:5) circle (.1);
\draw[fill] (306:5) circle (.1);
\draw[fill] (324:5) circle (.1);
\draw[fill] (342:5) circle (.1);
\draw(  -9:6) node {$\underline{2}5$};
\draw( -27:6) node {$\underline{1}24^{-1}$};
\draw( -45:6) node {$\underline{2}45^{-1}$};
\draw( -63:6) node {$\underline{1}4$};
\draw( -81:6) node {$1\underline{3}5$};
\draw( -99:6) node {$\underline{1}34^{-1}$};
\draw(-117:6) node {$\underline{3}5^{-1}$};
\draw(-135:6) node {$\underline{2}45$};
\draw(-153:6) node {$\underline{2}5^{-1}$};
\draw(-171:6) node {$\underline{3}5$};
\draw(-189:6) node {$\underline{2}35$};
\draw(-207:6) node {$\underline{2}4^{-1}$};
\draw(-225:6) node {$\underline{1}34$};
\draw(-243:6) node {$\underline{1}4^{-1}$};
\draw(-261:6) node {$\underline{2}4$};
\draw(-279:6) node {$\underline{1}24$};
\draw(-297:6) node {$\underline{1}3^{-1}$};
\draw(-315:6) node {$\underline{2}35^{-1}$};
\draw(-333:6) node {$\underline{1}3$};
\draw(-351:6) node {$1\underline{3}5^{-1}$};

\draw[thick] (  -9:5) to[bend right=30] (-153:5);
\draw[thick] ( -27:5) to[bend  left=20] (-279:5);
\draw[thick] ( -45:5) to[bend right=10] (-135:5);
\draw[thick] ( -63:5) to[bend right=15] (-243:5);
\draw[thick] ( -81:5) to[bend  left=20] (-351:5);
\draw[thick] ( -99:5) to[bend right=20] (-225:5);
\draw[thick] (-117:5) to[bend right=60] (-171:5);
\draw[thick] (-189:5) to[bend right=20] (-315:5);
\draw[thick] (-207:5) to[bend right=60] (-261:5);
\draw[thick] (-297:5) to[bend right=50] (-333:5);

\end{tikzpicture}
\caption{The relation $R_5$. Notation: e.g. $1\underline{3}5^{-1}$ denotes $\widehat{L}(3,\{1,3,5\})^{-1}$.}
\label{fig:R5_chord_diagram}
\end{figure}
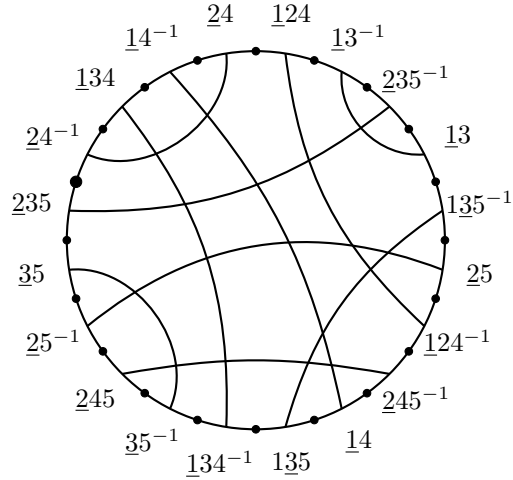
Note that this element belongs to the commutator subgroup as predicted by Proposition \ref{prp:relations_are_commutators}. The presentation complex of this relation is obtained from the $20$-gon by gluing the edges as in Figure \ref{fig:R5_chord_diagram} while preserving the orientation. This complex is easily verified to be a surface of genus $5$.

The Python implementation of an algorithm which computes the relation $R_m$ (and, more generally, the relation $\Red(R_\g(\lambda,J))\in F(\widehat{D})$ for any given $\K$, $J$ and $\lambda$) is available on the author's Github repository \cite{github}.

The complexity of the reduction grows exponentially on $m,$ but the algorithm is relatively fast: we are able to calculate the relation $R_m$ for $m=20$ in less then $30$ minutes. This relation is a word of length $8388612$.
\end{exm}
\begin{exm}
Consider an arbitrary presentation of the surface group $\pi_1(S_g)$ by $2g$ generators $x_1,\dots,x_{2g}$ modulo one relation $R$. Every generator appears in $R$ at least twice, since $R$ belongs to the commutator subgroup of $F(x_1,\dots,x_{2g})$ by Proposition \ref{prp:rel_in_comm}. Therefore, $|R|\geq 4g$. In the standard presentation $$\pi_1(S_g)\cong\langle a_1,b_1,\dots,a_g,b_g\mid a_1b_1a_1^{-1}b_1^{-1}\cdot\dotso\cdot a_gb_ga_g^{-1}b_g^{-1}\rangle$$  the relation satisfies $|R|=4g.$ Computer experiments show that, for the generators $\widehat{D}=\{\widehat{L}(i,J)\}$ of $\RC'_{\K_m}\cong \pi_1(S_{g(m)}),$ the length of the relation $R_m$ given by Theorem \ref{thm:explicit_presentation} is also equal to $4g(m)$ for $m\leq 20$ (and, possibly, for all $m$). Probably, this can be explained by the geometric nature of the proof of Theorem \ref{thm:explicit_presentation}: at the last step, we glue the disc $\ccK_{\underline{-1}}\simeq D^2$ to the punctured surface $(\R_{\K_m})^\circ:=\bigcup_{\g\neq\underline{-1}}\ccK_\g$ by the boundary, obtaining the closed surface $\R_{\K_m}$. We guess that the generators $\{\widehat{L}(i,J)\}$ of the free group $\pi_1((\R_{\K_m})^\circ)$ can be represented by simple closed curves in $\R_{\K_m}$ and hence $\R_{\K_m}$ is homeomorphic to the presentation complex for this presentation.

Our computer experiments show that the single relation between the Panov--Veryovkin generators $\{\widehat{\Gamma}(i,J)\}$ has length much larger than $4g(m).$ For example, if $m=10$ then $4g(m)=3076,$ but the relation between the Panov--Veryovkin generators has length $63940.$
\end{exm}

\section{Proof of Theorem \ref{thm:explicit_presentation}}
\label{section:upper_bound}

\subsection{Outline of the proof}
Consider the polyhedral product $B:=(\cone(\G),\G)^\K\simeq (E\G,\G)^\K$. In Proposition \ref{prp:loop_identification} we give a geometric interpretation of the isomorphism $\pi_1(B)\cong\Cart(\G,\K)$, identifying certain loops in $B$ with elements of $\Cart(\G,\K)$. Then, following Li Cai \cite{licai_slides}, we represent $B$ as the union of contractible subspaces $\cc(\K)_\g$ (see \eqref{eqn:ccKg_definition}) over all $\g\in\prod_{i=1}^m G_i$.

For a subset $Q\subset\prod_{i=1}^m G_i$, denote $B_Q:=\bigcup_{\h\in Q}\ccK_\h$. Then $B_{Q\sqcup\{\g\}}=B_Q\cup\ccK_\g,$ where $\ccK_\g$ is contractible. In some cases the intersection $B_Q\cap\ccK_\g$ is homotopy equivalent to $|\K_J|$ for $J=\supp\g\subset[m]$ (Proposition \ref{prp:admissible_intersection}). 
Hence, by an application of the van Kampen theorem (Lemma \ref{lmm:van_kampen_use}), the group $\pi_1(B_{Q\sqcup\{\g\}})$ is obtained from $\pi_1(B_{Q})$ by imposing $\rank\Pi_1(\K_J)$ relations and then adding $\b_0(\K_J)$ generators. By induction, in Theorem \ref{thm:explicit_presentation_admissible_subset} we give presentations for the fundamental group of $B_Q$ whenever $Q\subset\prod_{i=1}^mG_i$ is an \emph{admissible} subset (see Definition \ref{dfn:admissible}). A presentation of $\pi_1(B)$ is then obtained by passing to the limit (Lemma \ref{lmm:compactness_argument}).

\subsection{The geometric model for the classifying space}
\begin{figure}[h]
	\centering
	\begin{minipage}{0.4\textwidth}
		\centering
		\begin{tikzpicture}[scale=0.5, inner sep=2mm]
		
		\coordinate (1) at (-3,0);
		\coordinate (2) at (-3,1);
		\coordinate (3) at (-2,0);

		\fill (1) circle (2.5pt);
		\fill (2) circle (2.5pt);
		\fill (3) circle (2.5pt);

		\draw (1) -- (2);
		
		\draw (1) node[left] {$\ZZ_2$};
		\draw (1) node[below] {1};
		\draw (2) node[left] {$\ZZ_2$};
		\draw (2) node[above] {2};
		\draw (3) node[right] {$\ZZ_2$};
		\draw (3) node[below] {3};

		\draw[fill=gray!10] (0,0)--(1,5)--(6,5)--(5,0)--(0,0)--cycle;		
		\draw[fill=gray!10] (0,6)--(1,11)--(6,11)--(5,6)--(0,6)--cycle;
		\draw[dashed] (1,11)--(1,6);

		\draw (1,6)--(1,5);
		\draw (5,0)--(5,6);
		\draw (6,5)--(6,11);
		\draw (0,0)--(0,6);
		\fill (0,0) circle (4pt);
		\fill (1,5) circle (4pt);
		\fill (0,6) circle (4pt);
		\fill (1,11) circle (4pt);
		\fill (5,0) circle (4pt);
		\fill (6,5) circle (4pt);
		\fill (5,6) circle (4pt);
		
		\draw[ultra thick, fill=gray!40] (6,11)--(3.5,11)--(3,8.5)--(5.5,8.5)--(6,11)--cycle;
		\draw[ultra thick] (6,11)--(6,8);
		\fill (6,11) circle (6pt);

		\end{tikzpicture} 
	\end{minipage}
	\begin{minipage}{0.4\textwidth}
		\centering
		\begin{tikzpicture}[scale=0.5, inner sep=2mm]
		
		\coordinate (1) at (-3,0);
		\coordinate (2) at (-3,1);
		\coordinate (3) at (-2,0);

		\fill (1) circle (2.5pt);
		\fill (2) circle (2.5pt);
		\fill (3) circle (2.5pt);

		\draw (1) -- (2);
		
		\draw (1) node[left] {$\ZZ_2$};
		\draw (1) node[below] {1};
		\draw (2) node[left] {$\ZZ_3$};
		\draw (2) node[above] {2};
		\draw (3) node[right] {$\ZZ_2$};
		\draw (3) node[below] {3};

		\draw[fill=gray!30] (0,0)--(-1,2)--(4,2)--(5,0)--(0,0)--cycle;		
		\draw[fill=gray!30] (0,6)--(-1,8)--(4,8)--(5,6)--(0,6)--cycle;
		\draw[fill=gray!20] (-1,2)--(-2,3)--(3,3)--(4,2)--(-1,2)--cycle;
		\draw[fill=gray!20] (-1,8)--(-2,9)--(3,9)--(4,8)--(-1,8)--cycle;
\draw[fill=gray!10] (-0.33,3)--(1,5)--(6,5)--(4,2)--(3,3)--(-0.33,3)--cycle;
\draw[fill=gray!10] (-0.33,9)--(1,11)--(6,11)--(4,8)--(3,9)--(-0.33,9)--cycle;
		\draw[dashed] (-1,2)--(-0.33,3);
		\draw[dashed] (-1,8)--(-0.33,9);
		\draw[dashed] (1,11)--(1,6);
		\draw (1,6)--(1,5);
		\draw (3,3)--(3,9);
		\draw (-2,3)--(-2,9);
		\draw (5,0)--(5,6);
		\draw (6,5)--(6,11);
		\draw (0,0)--(0,6);

		\fill (0,0) circle (4pt);
		\fill (-2,3) circle (4pt);
		\fill (1,5) circle (4pt);
		\fill (0,6) circle (4pt);
		\fill (-2,9) circle (4pt);
		\fill (1,11) circle (4pt);
		\fill (5,0) circle (4pt);
		\fill (3,3) circle (4pt);
		\fill (6,5) circle (4pt);
		\fill (5,6) circle (4pt);
		\fill (3,9) circle (4pt);
		
		\draw[fill=gray!40] (6,11)--(3.5,11)--(2.1667,9)--(3,9)--(4,8)--(6,11)--cycle;
		\draw[fill=gray!70] (2.1667,9)--(1.5,8)--(4,8)--(3,9)--(2.1667,9)--cycle;
		\draw[ultra thick] (2.1667,9)--(3.5,11)--(6,11)--(4,8)--(1.5,8);
		\draw[ultra thick] (2.1667,9)--(1.5,8);
		\draw[ultra thick] (6,11)--(6,8);
		\draw (3,8)--(3,9);
		\fill (3,9) circle (4pt);
		\fill (6,11) circle (6pt);
		\end{tikzpicture} 
	\end{minipage}\qquad	
	\caption{Polyhedral products $B=(\cone\G,\G)^\K$ as unions of $\ccK_\g$}
	\label{fig:cell_models}
\end{figure}
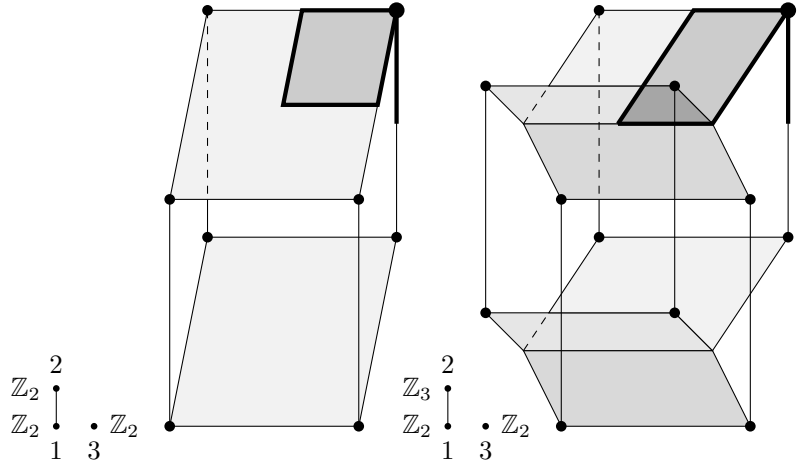

By Proposition \ref{prp:polyhedral_products_are_classifying_spaces}, the polyhedral product $(E\G,\G)^\K$ is the classifying space of the Cartesian subgroup $\Cart(\G,\K)$. Since $E\G$ is contractible, we have homotopy equivalences $(EG_i,G_i)\simeq (\cone G_i,G_i)$ and $(E\G,\G)^\K\simeq B$, where $B:=(\cone\G,\G)^\K=(\underline{X},\underline{A})^\K$ for $X_i:=\cone G_i$ and $A_i:=G_i$. Therefore, $(\cone\G,\G)^\K=B(\Cart(\G,\K)).$

We denote the vertex of $\cone G_i$ by $0$ and the segment between $0$ and $g_i$ by $[0,g_i]$. Therefore $\cone G_i=\bigcup_{g_i\in G_i}[0,g_i]$. For each $\g\in\prod_{i=1}^m G_i$, denote
\begin{equation}
\label{eqn:ccKg_definition}
\ccK_\g:=([0,\g],\{\g\})^\K.
\end{equation}
More formally, $\ccK_\g=(\underline{X},\underline{A})^\K$ for $X_i=[0,g_i],~A_i=\{g_i\}.$ It is clear that $B$ is glued out of $|\prod_{i=1}^m G_i|$ copies of $\ccK$,
$$(\cone\G,\G)^\K=\bigcup\Big\{\ccK_\g:\g\in\prod_{i=1}^m G_i\Big\}.$$
In particular, the real moment-angle complex $(\cone\ZZ_2,\ZZ_2)^\K=([-1,1],\{\pm 1\})^\K$ is glued out of $2^m$ copies of $\ccK.$

The polyhedral product $(\cone\G,\G)^\K$ and its subspace $\ccK_{\underline{1}}$ is shown in Figure \ref{fig:cell_models} for the case $\K=\bullet\!\!\!-\!\!\!\bullet\bullet$ and $\G=(\ZZ_2,\ZZ_2,\ZZ_2)$ (left) and $\G=(\ZZ_2,\ZZ_3,\ZZ_2)$ (right).

Throughout the proof, we will study the fundamental groups of path connected subspaces $B_Q:=\bigcup_{\g\in Q}\ccK_\g\subset (\cone\G,\G)^\K$ for certain $Q\subset\prod_{i=1}^mG_i$. Adding an element $\h$ to the set $Q$ corresponds to attaching a contractible space $\ccK_\h$ to the space $B_Q$ along a common CW-subcomplex. This affects the fundamental group as follows.

\begin{lmm}
\label{lmm:van_kampen_use}
Let $Y$ be a CW complex and $X,C\subset Y$ be its subcomplexes such that $X$ is path connected, $C$ is simply connected and $Y=X\cup C$. Let $X\cap C=A=\bigsqcup_{\alpha\in PC}A_\alpha$ be the decomposition into path components. Then a presentation of $\pi_1(Y)$ can be obtained from an arbitrary presentation of $\pi_1(X)$ by adding $\rank\Pi_1(A)$ relations and then adding $\b_0(A)$ generators, as follows.

Let $x\in X$ be a basepoint. For each $\alpha\in PC,$ choose a basepoint $a_\alpha\in A_\alpha$, a path $t_\alpha$ in $X$ from $x$ to $a_\alpha$, and a set of loops $\{\lambda_{i,\alpha}:i\in I_\alpha\}$ in $A_\alpha$ that generate the group $\pi_1(A_\alpha,a_\alpha)$, see Fig. \ref{fig:van_kampen_use}.

Choose a distinguished element $0\in PC.$ For each $\alpha\in PC\setminus\{0\},$ choose a path $s_\alpha$ in $C$ from $a_0$ to $a_\alpha.$

Then $\pi_1(Y,x)$ is the free product of the quotient group
$$\pi_1(X,x)/(t_\alpha\cdot \lambda_{i,\alpha}\cdot t_\alpha^{-1}:\alpha\in PC,~i\in I_\alpha)$$
with the free group generated by the loops $\{t_0\cdot s_\alpha\cdot t_\alpha^{-1}:\alpha\in PC\setminus\{0\}\}.$
\end{lmm}
\begin{figure}[h!]
\begin{tikzpicture}
  
	\draw[rounded corners=4mm, fill=gray!10] (0.1,0)--(1.5,0)--(1.5,1)--(2.4,1)--(2.4,0)--(4.9,0)--(4.9,2)--(0.1,2)--cycle;
  
	\coordinate (a0) at (0.7,0.5);  
	\coordinate (ai) at (3.2,0.5);
	\coordinate (x) at (2,-1.5);  
 	\coordinate (b) at (4.2,0.5);
 	\draw[rounded corners=4mm] (0,1)--(1,1)--(1,0)--(2.5,0)--(2.5,1)--(5,1)--(5,-2)--(0,-2)--cycle;
  
	\draw[fill] (a0) circle(0.05) node[left] {$a_0$};
	\draw[fill] (ai) circle(0.05) node[left] {$a_\alpha$};
	\draw[fill] (x) circle(0.05) node[above] {$x$};

	\draw[thick, directed] (x) to[out=180,in=270] (a0);
	\draw[thick, directed] (a0) to[out=90, in=90] (ai);
	\draw[thick, directed] (x) to[out=0, in=270] (ai);
	\draw[thick, directed] (ai) to[out=45, in=90] (b);
	\draw[thick] (b) to[out=-90,in=-45] (ai);
	
	\draw (1.0,-0.9) node[right] {\tiny $t_0$};
	\draw (3.4,-0.9) node[left] {\tiny $t_\alpha$};
	\draw (1.9,1.55) node[below] {\tiny $s_\alpha$};
	\draw (b) node[right] {\tiny $\lambda_{i,\alpha}$};
	\draw (0.1,1.7) node[right] {$C$};
	\draw (0,-1.7) node[right] {$X$};
	\draw (0,0.5) node[left] {$A_0$};
	\draw (5,0.5) node[right] {$A_\alpha$};

	\draw (3.5,0.55) .. controls (3.55,0.45) and (3.6,0.4) .. (3.7,0.4) .. controls (3.8,0.4) and (3.85,0.45) .. (3.9,0.55);

	\draw (3.525,0.5) .. controls (3.55,0.5) and (3.6,0.6) .. (3.7,.6) .. controls (3.8,0.6) and (3.85,0.5) .. (3.875,0.5);
\end{tikzpicture}
\caption{The space $Y=X\cup_A C$ from Lemma \ref{lmm:van_kampen_use}.}
\label{fig:van_kampen_use}
\end{figure}
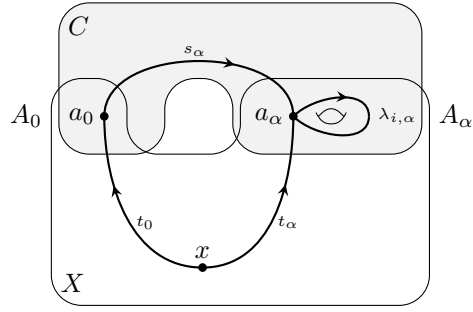
\begin{proof}
Without loss of generality, $x$ and $a_\alpha$ are 0-cells and $s_\alpha$ and $t_\alpha$ are distinct 1-cells of $Y.$ (If this is not the case, we first use cellular approximation theorem to ensure that images of $s_\alpha$ and $t_\alpha$ are in $\sk_1Y,$ and then attach a 1-cell and a 2-cell to each path without changing homotopy types, compare with the proof of \cite[Proposition 1.26]{hatcher}.) The closures of these cells are CW-subcomplexes $\overline{t_\alpha}=\{x,a_\alpha\}\cup t_\alpha\subset X$ and $\overline{s_\alpha}=\{a_0,a_\alpha\}\cup s_\alpha\subset C.$

Now denote $T:=\bigcup_{\alpha\in PC}\overline{t_\alpha}.$ We have the following pushout diagram of inclusion of connected CW-complexes with the basepoint $x$:
$$\xymatrix{
A\cup T\ar[r]^{i}
\ar[d]_-{j}
&
X
\ar[d]
\\
C\cup T
\ar[r]
&
Y.
}$$
Now we describe their fundamental groups and maps between them.
\begin{itemize}
\item Clearly, $(A\cup T)/T\cong\bigvee_{\alpha\in PC}A_\alpha.$ 
Since $T$ is contractible, we have
$$\pi_1(A\cup T,x)\cong\pi_1((A\cup T)/T)\cong\pi_1(\vee_{\alpha\in PC}A_\alpha)\cong\Bigast_{\alpha\in PC}\pi_1(A_\alpha,a_\alpha)=\Pi_1(A).$$
The isomorphism $\Pi_1(A)\overset\cong\longrightarrow\pi_1(A\cup T,x)$ is described as follows: if $\lambda$ is a loop in $A_\alpha$ with the basepoint $a_\alpha,$ then $t_\alpha\cdot\lambda\cdot t_\alpha^{-1}$ is the corresponding loop in $A\cup T$ with the basepoint $x.$ The composite map $i_*:\Pi_1(A)\cong\pi_1(A\cup T,x)\to \pi_1(X,x)$ has the same description.
\item We show that $\pi_1(C\cup T,x)$ is freely generated by the set of loops $\{t_0\cdot s_\alpha\cdot t_\alpha^{-1}:\alpha\in PC\setminus\{0\}\}$. Indeed, consider the contractible CW-subcomplex $T'=\bigcup_{\alpha\in PC\setminus\{0\}}\overline{s_\alpha}$ of $C$. Then $(C\cup T)/T'\cong C\vee\bigvee_{\alpha\in PC\setminus\{0\}}S^1$, where the circles are images of the loops $\{t_0\cdot s_\alpha\cdot t_\alpha^{-1}:\alpha\in PC\setminus\{0\}\}$ in $C\cup T$. Since $\pi_1(C)=1$, the group $\pi_1(C\cup T,x)\cong\pi_1((C\cup T)/T',x)$ is freely generated by these loops by the van Kampen theorem.
\end{itemize}
Applying the van Kampen theorem to the diagram above, we obtain the following pushout in the category of groups:
$$\xymatrix{
\Bigast_{\alpha\in PC}\pi_1(A,a_\alpha)
\ar[r]^-{i_*}
\ar[d]^-{j_*}
&
\pi_1(X,x)
\ar[d]
\\
F(\{t_0\cdot s_\alpha\cdot t_\alpha^{-1}:\alpha\in PC\setminus\{0\}\})
\ar[r]
&
\pi_1(Y,x),
}$$
where $i_*$ is the free product of maps $$\pi_1(A_\alpha,a_\alpha)\to \pi_1(X,x),\quad\lambda\mapsto t_\alpha\cdot \lambda\cdot t_\alpha^{-1},$$
while $j_*$ is the free product of composite maps
$$\pi_1(A_\alpha,a_\alpha)\to\underbrace{\pi_1(C,a_0)}_{=1}\to \pi_1(C\cup T,x),\quad \lambda\mapsto s_\alpha\cdot \lambda\cdot s_\alpha^{-1}\mapsto t_0\cdot s_\alpha\cdot\lambda\cdot s_\alpha^{-1}\cdot t_0^{-1}.$$
Since $\pi_1(C)=1$, the map $j_*$ is trivial.
Finally, each group $\pi_1(A_\alpha,a_\alpha)$ is generated by the set $\{\lambda_{i,\alpha}:i\in I_\alpha\}$, so the pushout is isomorphic to the group
\[
F(t_0\cdot s_\alpha\cdot t_\alpha^{-1}:\alpha\in PC\setminus\{0\})\,\ast\,\pi_1(X)/(t_\alpha\cdot\lambda_{i,\alpha}\cdot t_\alpha^{-1}:\alpha\in PC,i\in I_\alpha).\qedhere
\]
\end{proof}

\subsection{Subspaces of $\cc(\K)$}
The \emph{barycentric subdivision} of a simplicial complex $\K$ is the following simplicial complex $\K'$ on the vertex set $\K\setminus\{\varnothing\}$,
$$\K':=\Big\{\{I_1,\dots,I_r\}:~r\geq 0;~I_j\in\K,~j=1,\dots ,r;~I_1\subsetneq I_2\subsetneq\dots\subsetneq I_r;~I_1\neq\varnothing\Big\}.$$
Hence the faces of $\K'$ correspond to the chains $(I_1\subsetneq\dots\subsetneq I_r)$ of non-empty simplices in $\K$. The vertices of $\K'$ correspond to the non-empty simplices $I\in\K;$ we denote them by $(I)\in\K'$ and call the \emph{barycenters}. There is a natural PL homeomorphism $\beta:|\K'|\overset\cong\longrightarrow|\K|$ which maps $(I)\in\K'$ to the geometric barycenter $b_I$ of the geometric simplex $|I|=\conv(e_i:i\in I)\subset |\K|.$

The \emph{simplicial cone} over $\K'$ is naturally identified with the simplicial complex
$$\cone(\K'):=\{(I_1\subsetneq\dots\subsetneq I_r): r\geq 0;~I_j\in\K,~j=1,\dots,r\}$$
on the vertex set $\K.$ Note that $(\varnothing)$ is a vertex of $\cone(\K')$, and its geometric realisation is the apex $b_\varnothing$ of the geometric cone $|\cone(\K')|$ over $|\K'|$ (``the barycenter of the empty simplex'').

Now let $e_1,\dots,e_m$ be the standard basis of $\mathbb{R}^m,$ and denote $e_I:=\sum_{i\in I}e_i\in [0,1]^m$ for $I\subset[m].$

\begin{con}[{\cite[\S 2.9]{ToricTopology}}]
\label{con:cck}
Let $\K$ be a simplicial complex. Consider the PL embedding $\iota_\K:|\cone(\K')|\to[0,1]^m$ which is defined on vertices as $b_I\mapsto e_{[m]\setminus I}$ and extends linearly onto simplices of $|\cone(\K')|$.
\end{con}

\begin{lmm}[see {\cite[Proposition 2.9.12]{ToricTopology}}]
\label{lmm:cck_polyhedral}
The map $\iota_\K$ induces a homeomorphism $|\cone(\K')|\overset\cong\longrightarrow\ccK\subset[0,1]^m$, where $\ccK=([0,1],\{1\})^\K.$\qed
\end{lmm}

Recall that, for a simplex $A$ of a simplicial complex $\L,$ the subcomplex $$\st_\L A:=\{I\in\K:I\cup A\in\L\}$$ of $\L$ is called the \emph{star} of $J$. In the next two lemmas, we consider $\L=\K'$ and vertices $(\{j\})$ of $\K'$ that correspond to vertices $j$ of $\K$.
\begin{lmm}
For each $j\in[m]$, we have
$\iota_\K^{-1}(\{x_j=0\})=|\st_{\K'}(\{j\})|.$
\end{lmm}
\begin{proof}
Since $\iota_\K:|\cone(\K')|\to[0,1]^m$ is linear on simplices, the preimage of $H_j:=\{x\in[0,1]^m:x_j=0\}$ under $\iota_\K$ is (the geometric realisation of) the full subcomplex of $\cone(\K')$ spanned by the vertices that are mapped into $H_j$.

A vertex $b_I$ of $\cone(\K')$ is mapped to $e_{[m]\setminus I},$ hence $\iota_\K(b_I)\in H_j$ if and only if $j\in I$. It follows that $\iota_\K^{-1}(H_j)$ is the geometric realisation of the simplicial complex $\{(I_1\subsetneq\dots\subsetneq I_r)\in\cone(\K'):j\in I_1\},$ which clearly coincides with $\st_{\K'}(\{j\}).$
\end{proof}
\begin{lmm}\label{lmm:pl_things}
For $J\subset[m],$ denote $\mathcal{N}(J):=\bigcup_{j\in J}\st_{\K'}(\{j\})\subset\K'.$ Then
$$\cc(\K)\cap\bigcup_{j\in J}\{x_j=0\}=\iota_\K(|\mathcal{N}(J)|)\simeq|\K_J|.$$ 
\end{lmm}
\begin{proof}
We use $\mathcal{N}$ for $\mathcal{N}(J)$ throughout the proof. By the previous lemma, we have $\cc(\K)\cap\{x_j=0\}=\iota_\K(|\st_{\K'}(\{j\})|)$, so we obtain the first statement by taking the union over $j\in J$. 

Since $\iota_\K:|\cone(\K')|\to\ccK$ is a homeomorphism, now it is sufficient to prove that $|\mathcal{N}|$ and $|\K_J|$ are homotopy equivalent. Note that $\beta:|\K|\to|\K'|$ is a PL homeomorphism, and $\beta(|\K_J|)\subset|\mathcal{N}|$. 

Consider the full subcomplex of $\K'$ on the vertex set $\{(I):I\in\K,I\cap J\neq\varnothing\}$. Since $I\cap J\neq\varnothing$ if and only if $\{j\}\subset I$ for some $j\in J$, this full subcomplex coincides with $\mathcal{N}$.

Consider the full subcomplex of $\K'$ on the vertex set $\{(I):I\in\K,I\subset J\}$. This vertex set is equal to $\K_J\setminus\{\varnothing\},$ hence the geometric realisation of this full subcomplex coincides with $\beta(|\K_J|)$.

Hence we can define a strict deformation retraction of $|\mathcal{N}|$ onto $\beta(|\K_J|)$ as follows: for every $I\in\K$ such that $I\cap J\neq\varnothing,$ we move $b_I$ linearly to the point $b_{I\cap J}$ along the segment $|(I\cap J\subset I)|\subset|\mathcal{N}|.$  We extend this homotopy linearly onto geometric simplices of $|\mathcal{N}|$. This homotopy is well defined: if $(I_1\subsetneq\dots\subsetneq I_r)\in\mathcal{N}$, then the geometric simplex $\conv(b_{I_1},\dots,b_{I_r})\subset|\mathcal{N}|$ is linearly deformed onto the geometric simplex $\conv(b_{I_1\cap J},\dots,b_{I_r\cap J})$; the image of the homotopy lies in $\conv(b_{I_1},\dots,b_{I_r},b_{I_1\cap J},\dots,b_{I_r\cap J})\subset\beta(|I_r|)\subset|\mathcal{N}|.$
\end{proof}

\subsection{Admissible subsets}
For $\g\in\prod_{i=1}^mG_i$, denote $\supp\g:=\{i\in[m]:g_i\neq 1_i\}.$

For $\g\in\prod_{i=1}^mG_i$ and $I\subset[m]$, denote
$$\g(I):=(x_1,\dots,x_m)\in\prod_{i=1}^m G_i,~x_i=\begin{cases}
g_i,&i\in I;\\
1_i,&i\notin I.\end{cases}$$

\begin{dfn}
\label{dfn:admissible}
A subset $Q\subset\prod_{i=1}^m G_i$ is \emph{admissible} if $Q$ is finite, non-empty, and satisfies $\g(I)\in Q$ whenever $\g\in Q$ and $I\subset[m]$. If $Q$ is admissible, an element $\g\in Q$ is \emph{maximal} if we have $\g\neq\h(I)$ for all $\h\in Q\setminus\{\g\}$.
\end{dfn}

Consider the following partial ordering on $\prod_{i=1}^m G_i:$ $\g\leqslant\h$ if and only if $\g=\h(I)$ for some $I\subset[m]$. Clearly, $Q$ is admissible if and only if $Q$ is finite and is a \emph{lower set} (i.e. $\h\in Q$ whenever $\h\leqslant \g$ for some $\g\in Q$). Also, $\g\in Q$ is maximal if and only if $\g$ is a maximal element of $Q$ with respect to $\leqslant$. Hence the following properties are clear.
\begin{lmm}
\label{lmm:properties_of_admissible}
\begin{enumerate}
\item A finite union of admissible subsets is admissible.
\item For any $\g\in\prod_{i=1}^m G_i$, the subset $\{\g(I):I\subset[m]\}$ is admissible.
\item If $Q$ is admissible and $\g\in Q$ is maximal, then $Q\setminus\{\g\}$ is admissible.\qed
\end{enumerate}
\end{lmm}

For an admissible subset $Q\subset\prod_{i=1}^m G_i$, we denote $B_Q:=\bigcup_{\g\in Q}\ccK_\g$.

Denote the composition of homeomorphisms $\iota_\K:|\cone(\K')|\to\ccK$ and $\ccK\cong\ccK_\g$ by $\iota_\g$. In our notation, for $I\in\K$ we have
$$\iota_\g(b_I)=(x_1,\dots,x_m),\quad x_i=\begin{cases}
g_i,~i\notin I;\\
0,~i\in I.\end{cases}$$

\begin{lmm}
\label{lmm:ccKg_ccKh}
Let $x\in\ccK_\g\cap\ccK_\h$ for some $\g,\h\in\prod_{i=1}^mG_i$, and suppose that $x_j\neq 0$ for all $j\in\supp\g$. Then $\g=\h(S)$ for $S=\{i\in[m]:x_i\neq 0\}$.
\end{lmm}
\begin{proof}
If $i\in S$, then $x_i\neq 0$. Since $x\in\ccK_\g\cap\ccK_\h$, we have $x_i\in [0,g_i]\cap [0,h_i]$; it follows that $g_i=h_i$.

If $j\notin S$, then $x_j=0$. We have $j\notin\supp\g$ by the assumption, hence $g_j=1_j$.
\end{proof}
\begin{prp}
\label{prp:admissible_intersection}
Let $\g$ be a maximal element of an admissible subset $Q\subset\prod_{i=1}^m G_i$. Then $B_{Q\setminus\{\g\}}\cap\ccK_\g=\iota_\g(|\mathcal{N}(\supp\g)|)$. In particular, this space is homotopy equivalent to $\K_{\supp\g}$ by Lemma \ref{lmm:pl_things}.
\end{prp}
\begin{proof}
Denote $J:=\supp\g$ and $X:=B_{Q\setminus\{\g\}}.$
In view of Lemma \ref{lmm:pl_things}, it is sufficient to prove that $X\cap\ccK_\g=\bigcup_{j\in J}\{x_j=0\}\cap\ccK_\g.$

First let $j\in J$. Since $Q$ is admissible, we have $\g(J\setminus j)\in Q$, hence $\ccK_{\g(J\setminus j)}\subset X$. It is straightforward to show that $\ccK_{\g(J\setminus j)}\cap\ccK_\g=\{x_j=0\}\cap\ccK_\g$. Thus $X\cap\ccK_\g\supseteq\bigcup_{j\in J}\{x_j=0\}\cap\ccK_\g.$

On the other hand, let $x=(x_1,\dots,x_m)\in X\cap\ccK_\g,$ and suppose that $x_j\neq 0$ for all $j\in J$. We have $x\in\ccK_\h$ for some $\h\in Q\setminus\{\g\}$. By Lemma \ref{lmm:ccKg_ccKh}, $\g=\h(S)$ for some $S\subset[m]$. This contradicts the maximality of $\g$. Hence $x_j=0$ for some $j\in J$. We proved that $X\cap\ccK_\g\subseteq\bigcup_{j\in J}\{x_j=0\}\cap\ccK_\g.$
\end{proof}

\subsection{Loops defined by words}
\label{subsec:loop_identification}
Here we give a geometric interpretation of the isomorphism $\pi_1((\cone\G,\G)^\K,\underline{1})\cong\Cart(\G,\K)$ from Proposition \ref{prp:polyhedral_products_are_classifying_spaces}.

For each $k_i,g_i\in G_i$, consider the standard path in the topological space $\cone G_i$ from the point $k_i$ to the point $k_ig_i$ along two straight segments, $k_i\leadsto 0\leadsto k_ig_i$.

\begin{dfn}
\label{dfn:path_from_a_word}
Given a sequence of elements $g_{i_1},g_{i_2},\dots,g_{i_N}$ ($i_t\in[m]$, $g_{i_t}\in G_{i_t}$), we define the following path in the space $(\cone\G,\G)^\K$ which starts at the point $\underline{1}=(1_1,\dots,1_m)$:
$$(1_1,\dots,1_{i_1},\dots,1_{i_2},\dots)\leadsto (1_1,\dots,g_{i_1},\dots,1_{i_2},\dots)\leadsto(1_1,\dots,g_{i_1},\dots,g_{i_2},\dots)\leadsto \dots.$$
Here at each step $t=1,\dots,N$ we change the $i_t$-th coordinate from $k_{i_t}$ to $k_{i_t}g_{i_t}$ by following the standard path, while all the other coordinates are fixed. Hence the path ends at the point $(h_1,\dots,h_m)\in(\cone\G,\G)^\K$, where $h_i=\prod_{t:i_t=i}g_{i_t}\in G_i$.
\end{dfn}
\begin{prp}
\label{prp:loop_identification}
Let $g_{i_1},\dots,g_{i_N}$ be a sequence of letters such that $\prod_{t:i_t=i}g_{i_t}=1_i$ for each $i\in[m]$. Let $\lambda:[0,1]\to(\cone\G,\G)^\K$ be the corresponding path from $\underline{1}$ to $\underline{1}$ given in Definition \ref{dfn:path_from_a_word}. Then the isomorphism $\pi_1((\cone\G,\G)^\K,\underline{1})\to\Cart(\G,\K)$ from Proposition \ref{prp:polyhedral_products_are_classifying_spaces} maps $[\lambda]$ to $g_{i_1}\cdot\dotso\cdot g_{i_N}.$
\end{prp}
\begin{proof}
The element $g_{i_1}\cdot\dotso\cdot g_{i_N}\in G_1\ast\dots \ast G_m\cong\pi_1(BG_1\vee\dots\vee BG_m)$ is represented by the concatenation of loops which correspond to the elements $g_i\in G_i\cong\pi_1(BG_i)$. Since the isomorphism $\pi_1((B\G,\ast)^\K)\cong\G^\K$ from Proposition \ref{prp:polyhedral_products_are_classifying_spaces} is natural with respect to maps of simplicial complexes, the word $g_{i_1}\cdot\dotso\cdot g_{i_N}\in\G^\K\cong\pi_1((B\G,\ast)^\K)$ is represented by the same loop. On the other hand, this loop is the composition of $\lambda:[0,1]\to (\cone\G,\G)^\K$ with the homotopy equivalence $(\cone\G,\G)^\K\simeq (E\G,\G)^\K$ and the natural map $(E\G,\G)^\K\to (B\G,\ast)^\K$. This map is injective on the fundamental groups, so the claim follows.
\end{proof}

\begin{rmk}
In fact, the construction above provides a bijective correspondence
$$\{x\in\G^\K:\pi(x)=\h\}\overset\cong\longrightarrow\{\text{paths in }(\cone\G,\G)^\K\text{ from }\underline{1}\text{ to }\h\}/\sim,~\forall\h\in\prod_{i=1}^m G_i,$$
where $\pi:\G^\K\to\prod_{i=1}^m G_i$ is the natural projection and $\sim$ is the homotopy modulo endpoints. Different factorizations $x=g_{i_1}\cdot\dotso\cdot g_{i_N}\in\G^\K$ correspond to homotopic paths: for example, each edge $\{i,j\}\in\K$ gives rise to the relation $g_ig_j=g_jg_i$ in $\G^\K$ (hence to ambiguity in the factorization) and to disjoint union of squares in $(\cone\G,\G)^\K$ (hence to homotopies between the corresponding paths). We do not prove this fact since we do not use it in full generality.
\end{rmk}

\begin{dfn}
To each standard element $L_\g(i,J)\in D$ we assign the loop in the space $(\cone\G,\G)^\K$ as in Definition \ref{dfn:path_from_a_word} using the sequence from Definition \ref{dfn:LiJ} (e.g. $L_\g(3,\{1,3,4\})$ is considered as the sequence $g_1g_3g_4g_3^{-1}g_4^{-1}g_1^{-1}$). Hence each word $w=L_\g(i_1,J_1)^{\pm 1}\cdot\dots\cdot L_\g(i_N,J_N)^{\pm 1}$ on the standard elements defines a loop in $(\cone\G,\G)^\K$. By Proposition \ref{prp:loop_identification}, this loop represents the element $w\in\Cart(\G,\K)$ under the isomorphism $\pi_1((\cone\G,\G)^\K,\underline{1})\cong\Cart(\G,\K)$.

We say that such a word $w$ \emph{defines} a loop in a subspace $X\subset(\cone\G,\G)^\K$ if the corresponding loop in $(\cone\G,\G)^\K$ lies in $X$. Hence if $L_\g(i,J)$ defines a loop in $X$, then it represents an element in $\pi_1(X,\underline{1})$, which is mapped into $L_\g(i,J)$ under the homomorphism $\pi_1(X,\underline{1})\to\pi_1((\cone\G,\G)^\K,\underline{1})\cong\Cart(\G,\K)$.
\end{dfn}

\subsection{Fundamental groups of subspaces corresponding to admissible subsets}
\begin{thm}
\label{thm:explicit_presentation_admissible_subset}
Let $Q\subset\prod_{i=1}^m G_i$ be an admissible subset. For each $\g\in Q$, choose a set of cycles $\Gen(\g)$ in $\K_{\supp\g}$ that generate the group $\Pi_1(\K_{\supp\g}).$

Then the group $\pi_1(B_Q,\underline{1})$ is presented by $\sum_{\g\in Q}\b_0(\K_{\supp\g})$ generators, which correspond to the words
$$\{\widehat{L}_\g(i,\supp\g):~\g\in Q,~i\in\Theta(\supp\g)\},$$
modulo $\sum_{\g\in Q}|\Gen(\g)|$ relations
$$\{\Red(R_\g(\lambda,\supp\g))=1:\g\in Q,~\lambda\in\Gen(\g)\}.$$
\end{thm}
($R_\g(\lambda,\supp\g)\in F(D)$ and $\Red:F(D)\to F(\widehat{D})$ were defined in Subsection \ref{subsection:statement_main}.)
\begin{rmk}
In particular, we claim that the word $\Red(R_\g(\lambda,\supp\g))\in F(\widehat{D})$ depends only on the generators $\{\widehat{L}_\h(i,\supp\h):\h\in Q,i\in\Theta(\supp\h)\}\subset\widehat{D}.$
\end{rmk}

\begin{proof}[Proof of Theorem \ref{thm:explicit_presentation_admissible_subset}]
Induction on $|Q|$. The base case is $Q=\{\underline{1}\}$, when the presentation is empty and $B_Q=\ccK_{\underline{1}}$ is contractible.

On the induction step, we choose a maximal element $\g\in Q$, and denote $X:=B_{Q\setminus\{\g\}}$, $J:=\supp\g$. Hence $B_Q=X\cup\ccK_\g,$ and $X\cap\ccK_\g\simeq|\K_J|$ by Proposition \ref{prp:admissible_intersection}. It is sufficient to prove that a presentation of $\pi_1(B_Q,\underline{1})$ can be obtained from a presentation of $\pi_1(X,\underline{1})$ by adding the $\b_0(\K_J)$ generators $\{\widehat{L}_\g(i,J):i\in\Theta(J)\}$ and the $|\Gen(\g)|$ relations $\{\Red(R_\g(\lambda,J)):\lambda\in\Gen(\g)\}$.

We first define some paths in $B_Q$ and show that, up to homotopy, their concatenations are defined by words on elements from $\widehat{D}\subset\Cart(\G,\K)$.
\begin{itemize}
\item For each $I\subset[m]$ consider the point $v_I:=\g([m]\setminus I)\in B_Q$ with coordinates
$(v_I)_j=\begin{cases}1_j,&j\in I\\g_j,&j\notin I\end{cases}$
and the path
$$p(I):(1_1,\dots,1_m)\leadsto ((v_I)_1,1_2,\dots,1_m)\leadsto\dots\leadsto ((v_I)_1,\dots,(v_I)_m)$$
in $B_Q$ from $\underline{1}$ to $v_I$. (Both $v_I$ and $p(I)$ are in $B_Q$, since $Q$ is admissible and $\g\in Q$. Moreover, for $I\neq\varnothing$ both $v_I$ and $p(I)$ are in $X$.)

\item For each $i\in[m]$ consider the point $a_i\in\ccK_\g$ with the coordinates
$(a_i)_k:=\begin{cases}
0,&k=i\\
g_k,&k\neq i
\end{cases}$
and the path $$t(i):
\underbrace{(1_1,\dots,1_m)}_{=\underline{1}}\overset{p(\{i\})}\leadsto\underbrace{(g_1,\dots,1_i,\dots,g_m)}_{=v_{\{i\}}}\leadsto\underbrace{(g_1,\dots,0_i,\dots,g_m)}_{=a_i}$$
in $X$ from $\underline{1}$ to $a_i$.
(Here the second part is the standard path along the segment $[1_i,0_i]$.)
\item For each $i\in J$ consider the path $s_i:a_{\max(J)}\leadsto v_\varnothing\leadsto a_i$ in $\ccK_\g$ from $a_{\max(J)}$ to $a_i$.  (Here $s_i$ the concatenation of the standard paths along the segments $[0_{\max(J)},g_{\max(J)}]$ and $[g_i,0_i]$).
\end{itemize}
\begin{lmm}
\label{lmm:identify_generators}
For $i\in J\subset[m]$, denote $m'=\max(J).$ Then the loop $t(m')\cdot s_i\cdot t(i)^{-1}:\underline{1}\leadsto a_{m'}\leadsto a_i\leadsto \underline{1}$ defines the element in $\pi_1(B_Q,\underline{1})$ represented by the word $L_\g(i,J)$.
\end{lmm}
\begin{proof}
By definition, this loop is the following composition of paths:
$$
\underline{1}=(1_1,\dots,1_i,\dots,1_{m'},\dots,1_m)\overset{p(\{m'\})}\leadsto
(g_1,\dots,g_i,\dots,1_{m'},\dots,g_m)\leadsto$$ 
$$(g_1,\dots,g_i,\dots,0_{m'},\dots,g_m)\leadsto (g_1,\dots,g_i,\dots,g_{m'},\dots,g_m)\leadsto$$
$$
(g_1,\dots,0_i,\dots,g_{m'},\dots,g_m)\leadsto
(g_1,\dots,1_i,\dots,g_{m'},\dots,g_m)\leadsto\dots\leadsto \underline{1}.$$
Since $g_j=1_j$ for $j>m'$, this path is defined by the sequence of letters $g_1g_2\dots g_{m'}\cdot g_i^{-1}\cdot g_{m'}^{-1}\dots g_{i+1}^{-1}g_{i-1}^{-1}\dots g_1^{-1}=L_\g(i,J).$
\end{proof}
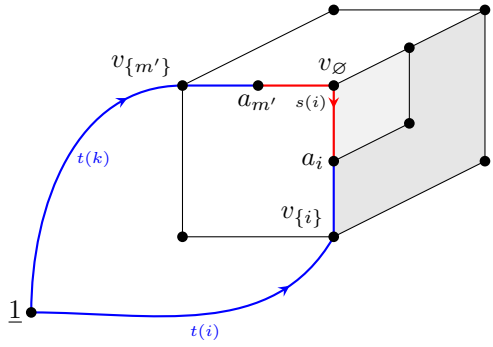
\begin{figure}[h!]
	\begin{tikzpicture}[scale=1]

	\coordinate (1) at (-2,-1);
	\coordinate (2) at (0,0);
	\coordinate (3) at (2,0);
	\coordinate (4) at (0,2);
	\coordinate (5) at (2,2);
	\coordinate (6) at (1,1);
	\coordinate (7) at (2,1);
	\coordinate (8) at (1,2);
	\coordinate (9) at (3,2.5);
	\coordinate (10) at (4,3);
	\coordinate (11) at (2,3);
	\coordinate (12) at (4,1);
	\coordinate (13) at (3,1.5);
	
	\fill[gray!20] (3)--(12)--(10)--(9)--(13)--(7)--cycle;
	\fill[gray!10] (13)--(9)--(5)--(7)--cycle;

	\draw (9)--(13)--(7)--(5)--(8);
	\draw (3)--(2)--(4);
	\draw (4)--(11)--(10)--(12)--(3);
	\draw (5)--(9)--(10);

	\draw[thick, blue, directed] (1) to[out=0,in=240] node[below] {\tiny $t(i)$}
	(3)--(7);
	
	\draw[thick, blue, directed] (1) to[out=90, in=180] node[right] {\tiny $t(k)$} (4) -- (8);
	
	\draw[thick, red, directed] (8)--(5)--(7);
	
	\foreach \i in {1, ..., 5} {\fill (\i) circle (2pt);}
	\foreach \i in {7, ..., 13} {\fill (\i) circle (2pt);}

	\draw (1) node[left] {$\underline{1}$};
	\draw (3) node[above left] {$v_{\{i\}}$};
	\draw (4) node[above left] {$v_{\{m'\}}$};
	\draw (5) node[above] {$v_{\varnothing}$};
	\draw (7) node[left] {$a_i$};
	\draw (8) node[below] {$a_{m'}$};

	\draw (1.7,2) node[below] {\tiny $s(i)$};
	
	\end{tikzpicture}
	\caption{Proof of Lemma \ref{lmm:identify_generators}}
	\label{fig:new_generator}
\end{figure}
\begin{itemize}
\item For each edge $\{i,j\}\in\K_J$, consider the point $u_{ij}\in\ccK_\g$ with coordinates
$(u_{ij})_k:=\begin{cases}
g_k,&k\neq i,j\\
0,&k\in\{i,j\}
\end{cases}
$
and the path $\ell(i,j):a_i\leadsto u_{ij}\leadsto a_j$ in $\ccK_\g$. 
\end{itemize}
\begin{lmm}
\label{lmm:identify_edge}
Let $\{i,j\}\in\K_J$. Then the loop $t(i)\cdot\ell(i,j)\cdot t(j)^{-1}:\underline{1}\leadsto a_i\leadsto a_j\leadsto \underline{1}$ defines the element in $\pi_1(X,\underline{1})$ represented by the word $L_\g(j,J\setminus i)\cdot L_\g(i,J\setminus j)^{-1}$.
\end{lmm}
\begin{proof}
By definition, $t(i)\cdot\ell(i,j)\cdot t(j)^{-1}$ is the concatenation
$\underline{1}\overset{p(\{i\})}\leadsto v_{\{i\}}\leadsto a_i\leadsto u_{ij}\leadsto a_j\leadsto v_{\{j\}}\overset{p(\{j\})^{-1}}\leadsto \underline{1}.$ The paths $v_{\{i\}}\leadsto a_i\leadsto u_{ij}\leadsto a_j\leadsto v_{\{j\}}$ and $v_{\{i\}}\leadsto v_{\{i,j\}}\leadsto v_{\{j\}}$ are homotopic in $X$, since there are three squares
$$\prod_{k\neq i,j}\{g_k\}\times ([0_i,1_i]\times[0_j,1_j]\cup[0_i,1_i]\times[0_j,g_j]\cup[0_i,g_j]\times[0_j,1_j])\subset X.$$
By adding and removing the path $p(\{i,j\})^{-1}\cdot p(\{i,j\})$ from $v_{\{i,j\}}$ to $\underline{1}$ and back, we see that our loop is homotopic in $X$ to the loop
$$\underline{1}\overset{p(\{i\})}\leadsto v_{\{i\}}\leadsto v_{\{i,j\}}\overset{p(\{i,j\})^{-1}}\leadsto\underline{1}\overset{p(\{i,j\})}\leadsto v_{\{i,j\}}\leadsto v_{\{j\}}\overset{p(\{j\})}\leadsto \underline{1}$$
which is represented by the word $$\prod_{k\neq i}g_i\cdot g_j^{-1}\cdot (\prod_{k\neq i,j}g_k)^{-1}\cdot \prod_{k\neq i,j}g_k\cdot g_i\cdot(\prod_{k\neq j}g_j)^{-1}=L_\g(J\setminus i,j)\cdot L_\g(J\setminus j,i)^{-1}.\qedhere$$
\end{proof}
\begin{figure}[h!]
	\centering
	\begin{minipage}{0.4\textwidth}
	\centering
	\begin{tikzpicture}[scale=1, inner sep=1mm]

	\coordinate (1) at (-2,-1);
	\coordinate (2) at (0,0);
	\coordinate (3) at (2,0);
	\coordinate (4) at (0,2);
	\coordinate (5) at (2,2);
	\coordinate (6) at (1,1);
	\coordinate (7) at (2,1);
	\coordinate (8) at (1,2);

	\draw[fill=gray!20] (2)--(3)--(7)--(6)--(8)--(4)--cycle;	
	\draw[dashed] (7)--(5)--(8);
	\draw[dashed] (0,1) -- (6) -- (1,0);
			
	\draw[thick, blue, directed] (1) to[out=0,in=240] (3);
	\draw[thick, blue, directed] (3) --
	(7) -- (6) -- (8) -- (4);
	\draw[thick, blue, directed] (4) to[out=180, in=90] (1);
	
	\foreach \i in {1, ..., 8} {\fill (\i) circle (2pt);}

	\draw (1) node[left] {$\underline{1}$};
	\draw (2) node[left] {$v_{\{i,j\}}$};
	\draw (3) node[right] {$v_{\{i\}}$};
	\draw (4) node[above] {$v_{\{j\}}$};
	\draw (5) node[above] {$v_{\varnothing}$};
	\draw (1,0.8) node[left] {$u_{i,j}$};
	\draw (7) node[right] {$a_i$};
	\draw (8) node[above] {$a_j$};

	\draw (0.1,-0.1) node[below] {\tiny $(1_i,1_j)$};
	\draw (2, 1.9) node[right] {\tiny $(g_i,g_j)$};
	\draw (0.95, 0.8) node[right] {\tiny $(0_i,0_j)$};

	\draw[blue] (0, -1) node[above] {\tiny $p(\{i\})$};
	\draw[blue] (-0.6, 2) node[left] {\tiny $p^{-1}(\{j\})$};
	
	\end{tikzpicture}
	\end{minipage}
	\begin{minipage}{0.4\textwidth}
	\centering
	\begin{tikzpicture}[scale=1, inner sep=1mm]

	\coordinate (1) at (-2,-1);
	\coordinate (2) at (0,0);
	\coordinate (3) at (2,0);
	\coordinate (4) at (0,2);
	\coordinate (5) at (2,2);
	\coordinate (6) at (1,1);
	\coordinate (7) at (2,1);
	\coordinate (8) at (1,2);

	\draw[fill=gray!20] (2)--(3)--(7)--(6)--(8)--(4)--cycle;	
	\draw[dashed] (7)--(5)--(8);
	\draw[dashed] (0,1) -- (6) -- (1,0);
			
	\draw[thick, blue, directed] (1) to[out=0,in=240] (3);
	\draw[thick, blue, directed] (3) -- (0.2,0);
	\draw[thick, blue, directed] (0.2,0) -- (-1.9,-0.95);
	\draw[thick, blue, directed] (-1.9,-0.95) -- (0,0.1);
	\draw[thick, blue, directed] (0,0.1) -- (4);

	\draw[thick, blue, directed] (4) to[out=180, in=90] (1);
	
	\foreach \i in {1, ..., 8} {\fill (\i) circle (2pt);}

	\draw (1) node[left] {$\underline{1}$};
	\draw (0,0.2) node[left] {$v_{\{i,j\}}$};
	\draw (3) node[right] {$v_{\{i\}}$};
	\draw (4) node[above] {$v_{\{j\}}$};
	\draw (5) node[above] {$v_{\varnothing}$};
	\draw (1,0.8) node[left] {$u_{i,j}$};
	\draw (7) node[right] {$a_i$};
	\draw (8) node[above] {$a_j$};
	\draw[blue] (1,0.05) node[below] {\tiny $g_j^{-1}$};
	\draw[blue] (-1,-0.55) node[right] {\tiny $p^{-1}(\{i,j\})$};

	\end{tikzpicture} 
	\end{minipage}
	\caption{Proof of Lemma \ref{lmm:identify_edge}}
	\label{fig:heartbreak}
\end{figure}
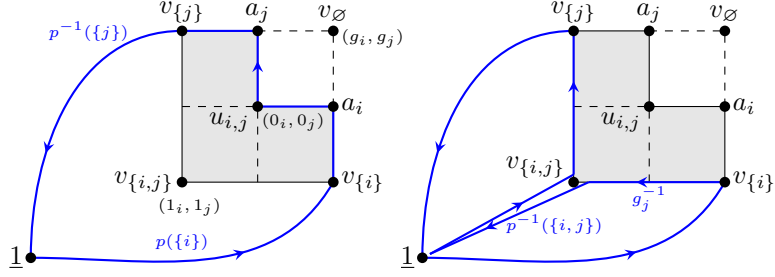
\begin{itemize}
\item For each simplicial cycle $\lambda=(i_1,\dots,i_k,i_{k+1}=i_1)$ in $\K_J$, consider the path $\ell(\lambda):a_{i_1}\overset{\ell(i_1,i_2)}\leadsto a_{i_2}\leadsto\dots\leadsto a_{i_{k}}\overset{\ell(i_k,i_1)}\leadsto a_{i_1}$ in $\ccK_\g$.
\end{itemize}
\begin{lmm}
\label{lmm:identify_relations_preliminary}
Let $\lambda=(i_1,\dots,i_k,i_{k+1}=i_1)$ be a cycle in $\K_J$. Then the loop
$$\gamma_\lambda=t(i_1)\cdot\ell(\lambda)\cdot t(i_1)^{-1}:\underline{1}\leadsto a_{i_1}\leadsto a_{i_1}\leadsto \underline{1}$$
defines the element in $\pi_1(X,\underline{1})$ represented by the word
$$R_\g(\lambda,J)=\prod_{t=1}^kL_\g(i_{t+1},J\setminus i_t)\cdot L_\g(i_t,J\setminus i_{t+1})^{-1}.$$
\end{lmm}
\begin{proof}
The loop $\gamma_\lambda$ is the concatenation of loops $t(i_1)\ell(i_1,i_2)\ell(i_2,i_3)\dots\ell(i_k,i_1)t(i_1)^{-1}$. By adding the trivial paths $t(i_s)^{-1}t(i_s)$ from $a_{i_s}$ to $\underline{1}$ and back, we see that $\gamma_\lambda$ is homotopic to the loop $\gamma'_\lambda=t(i_1)\ell(i_1,i_2)t(i_2)^{-1}\cdot t(i_2)\ell(i_2,i_3)t(i_3)^{-1}\cdot\dotso\cdot t(i_k)\ell(i_k,i_1)t(i_1)^{-1}.$ This is the concatenation of loops identified in Lemma \ref{lmm:identify_edge}.
\end{proof}
\begin{lmm}
\label{lmm:identify_relations_final}
Let $Q\subset\prod_{i=1}^mG_i$ be an admissible subset, and let $\g\in Q$, $J=\supp\g,$ $i\in J$. Then the words $L_\g(i,J)$ and $\Red(L_\g(i,J))$ define the same element in $\pi_1(B_{Q},\underline{1}).$
\end{lmm}
\begin{proof} We follow the algorithm from Proposition \ref{prp:reduction_of_generators} which computes $\Red(L_\g(i,J))$.
On each step we have a concatenation of loops in $\pi_1(B_Q,\underline{1})$, one of which corresponds to the generator $L_\g(i,J)$, and we replace it with the concatenation of loops $L_\g(j,J)\cdot L_\g(i,J\setminus j)\cdot L_\g(j,J\setminus i)^{-1}$ for some $\{i,j\}\in \K_{J}.$ It is sufficient to show that these two loops are equal in $\pi_1(B_Q,\underline{1})$; equivalently, that the relation \eqref{eqn:elementary_relation} holds for the corresponding paths in $\pi_1(B_Q,\underline{1})$. The left hand side of this relation is represented by the path $\underline{1}\leadsto v_\varnothing\leadsto v_{\{i\}}\leadsto \underline{1}\leadsto v_{\{i\}}\leadsto v_{\{i,j\}}\leadsto \underline{1}.$
This loop is homotopic in $B_Q$ to the loop
$\underline{1}\leadsto v_\varnothing\leadsto v_{\{i\}}\leadsto v_{\{i,j\}}\leadsto \underline{1},$ and the paths $v_\varnothing\leadsto v_{\{i\}}\leadsto v_{\{i,j\}}$ and $v_\varnothing \leadsto v_{\{j\}}\leadsto v_{\{i,j\}}$ are homotopic since $B_Q$ contains the four squares
$$\prod_{k\neq i,j}\{g,k\}\times ([0_i,1_i]\cup[0_i,g_i])\times([0_j,1_j]\cup[0_j,g_j]).$$
Repeating the argument, we obtain the right hand side of this relation.
\end{proof}
\begin{crl}
\label{crl:identification_of_loops}
The loop $t(\max(J))\cdot s_i\cdot t(i)^{-1}\in\pi_1(B_Q,\underline{1})$ represents the word $L_\g(i,J)$ whenever $i\in J$.

The loop $t(i_1)\ell(i_1,i_2)\dots\ell(i_k,i_1)t(i_1)^{-1}\in\pi_1(X,\underline{1})$ represents the word $\Red(R_\g(\lambda,J))$ whenever $\lambda=(i_1,\dots,i_k,i_1)$ is a simplicial cycle in $\K_J$.\qed
\end{crl}

We have $B_Q=X\cup\ccK_\g$, where $\ccK_\g$ is contractible and $X\cap\ccK_\g=:C=\iota_\g(|\mathcal{N}(J)|)\simeq|\K_J|$ by Proposition \ref{prp:admissible_intersection}. The path components of $|\K_J|$ are represented by the set $\Theta(J)\sqcup\{\max(J)\}$, where the last path component is distinguished; under the homotopy equivalence $X\cap\ccK_\g\simeq|\K_J|$, vertices $i\in J$ correspond to points $a_i\in\ccK_\g$, while edges $\{i,j\}\in\K_J$ correspond to paths $a_i\leadsto u_{ij}\leadsto a_j$ in $\ccK_\g$. We denote $m':=\max(J)$, and consider the following paths:

\begin{itemize}
\item For each $i\in\Theta(J)\sqcup\{m'\}$, the path $t_i:\underline{1}\overset{t(i)}\leadsto a_i$ in $X$;
\item For each $i\in\Theta(J)$, the path $s_i:a_{m'}\leadsto v_\varnothing\leadsto a_i$ in $\ccK_\g$;
\item For each cycle $\lambda=(i_1,\dots,i_k,i_{k+1}=i_1)$ in $\K_J$, the path $\ell(\lambda):a_{i_1}\overset{\ell(i_1,i_2)}\leadsto a_{i_2}\leadsto\dots\leadsto a_{i_{k}}\overset{\ell(i_k,i_1)}\leadsto a_{i_1}$ in $\ccK_\g$.
\end{itemize}

This set of paths satisfies the conditions of Lemma \ref{lmm:van_kampen_use} (here $PC=\Theta(J)\sqcup\{m'\}$ and $x=\underline{1}$). Then Corollary \ref{crl:identification_of_loops} expresses the resulting generators and relations in terms of words in $\Cart(\G,\K)$. This completes the proof of Theorem \ref{thm:explicit_presentation_admissible_subset}.
\end{proof}

\subsection{Infinite case}
Recall that a collection of spaces $\{X_\alpha\}_{\alpha\in I}$ is \emph{directed} if for any $\alpha_1,\alpha_2\in I$ there exists $\alpha\in I$ such that $X_{\alpha_1}\cup X_{\alpha_2}\subset X_\alpha.$ In the next lemma, for a presentation $G=\langle \Gamma\mid R\rangle$ we have $\Gamma\subset G$ and $R\subset F(\Gamma)$.

\begin{lmm}
\label{lmm:compactness_argument} Let $(X,x_0)$ be a pointed CW-complex and $\{X_\alpha\}$ be a directed set of CW-subcomplexes $X_\alpha\subset X$ such that $X=\bigcup_\alpha X_\alpha,$ $x_0\in\bigcap_\alpha X_\alpha$.

Suppose that there are presentations $\pi_1(X_\alpha,x_0)=\langle\Gamma_\alpha\mid R_\alpha\rangle$ such that for each inclusion $X_\alpha\subset X_\beta$ the homomorphism $i_{\alpha\beta}:\pi_1(X_\alpha,x_0)\to\pi_1(X_\beta,x_0)$ induces inclusions of presentations. (In other words, $i_{\alpha\beta}$ restricts to an inclusion $\Gamma_\alpha\hookrightarrow\Gamma_\beta$, and the induced map $(i_{\alpha\beta})_*:F(\Gamma_\alpha)\to F(\Gamma_\beta)$ restricts to an inclusion $R_\alpha\hookrightarrow R_\beta$.)

Then $\pi_1(X)\cong\Big\langle\bigcup_\alpha\Gamma_\alpha\Big|\bigcup_\alpha R_\alpha\Big\rangle.$
\end{lmm}
\begin{proof}
By \cite[Propositions A.4, A.5]{hatcher}, every point $x\in X$ has a contractible open neighbourhood $N(x)\subset X$ such that each CW-subcomplex $X_\alpha$ is a strict deformation retract of its neighbourhood $V_\alpha:=\bigcup_{x\in X_\alpha}N(x).$ Replacing $X_\alpha$ with $V_\alpha$, we now assume that $X_\alpha$ are open instead of being CW-subcomplexes.

Moreover, every compact subset $K\subset X$ belongs to some $X_{\alpha(K)}.$ Indeed, we have an open covering $K\subset\bigcup_\alpha X_\alpha;$ choose a finite subcovering $K\subset X_{\alpha_1}\cup\dots\cup X_{\alpha_n}.$ Since $\{X_\alpha\}$ is directed, we have $X_{\alpha_1},\dots,X_{\alpha_n}\subset X_\alpha$ for some $\alpha=\alpha(K).$

Now let $\gamma:(S^1,\pt)\to (X,x_0)$ represent an element $[\gamma]\in\pi_1(X,x_0).$ Then $\Img\gamma\subset X$ is compact, so $\Img\gamma\subset X_\alpha$ for some $\alpha,$ hence $[\gamma]$ belongs to the subgroup of $\pi_1(X)$ generated by $\Gamma_\alpha.$ It follows that $\pi_1(X)$ is generated by the set $\Gamma:=\bigcup_{\alpha}\Gamma_\alpha\subset\pi_1(X),$ i.e. the natural map $F(\Gamma)\to\pi_1(X)$ is surjective.

Similarly, let $r\in \Ker(F(\Gamma)\to\pi_1(X))$ be a relation. It depends only on a finite set of generators $[\gamma_1],\dots,[\gamma_N]\in\Gamma.$ For every $j=1,\dots,N$ we have $[\gamma_j]\in\Gamma_{\alpha_j}\subset\pi_1(X_{\alpha_j})$ for some $\alpha_j.$
The relation $r$ corresponds to a homotopy between a loop in $X$ and the trivial map. The image of homotopy is an image of unit square, so it is compact and is contained in some $X_{\alpha_0}.$ For some $\alpha$, we have $X_{\alpha_0},X_{\alpha_1},\dots,X_{\alpha_N}\subset X_\alpha.$ Thus the word $r$ is a relation between the elements $[\gamma_1],\dots,[\gamma_N]$ of $\pi_1(X_\alpha),$ hence it belongs to the normal closure of the subset $R_\alpha\subset F(\Gamma_\alpha).$ It follows that $r$ belongs to the normal closure of $\bigcup_\alpha R_\alpha\subset F(\bigcup_\alpha\Gamma_\alpha).$
\end{proof}

\begin{proof}[Proof of Theorem \ref{thm:explicit_presentation}]

If all groups $G_i$ are finite, it suffices to apply Theorem \ref{thm:explicit_presentation_admissible_subset} to the admissible subset $Q:=\prod_{i=1}^m G_i.$  In the general case, we apply Lemma \ref{lmm:compactness_argument} to the directed covering $\{B_Q:Q\subset\prod_{i=1}^mG_i\text{ is admissible}\}$ of $(B,\underline{1})$ and the presentations of $\pi_1(B_Q,\underline{1})$ provided by Theorem \ref{thm:explicit_presentation_admissible_subset}. (This is a directed covering by Lemma \ref{lmm:properties_of_admissible}.) Since any $\g\in\prod_{i=1}^mG_i$ belongs to an admissible subset, we obtain the set $$\{\widehat{L}_\g(i,\supp\g):\g\in\prod_{i=1}^mG_i,~i\in\Theta(\supp\g)\}=\{\widehat{L}_\g(i,J):J\subset[m],~\g\in\prod_{j\in J}G_j^*,~i\in\Theta(J)\}$$
 of generators for $\Cart(\G,\K)$, and similarly with the relations.
\end{proof}

\section{Discussion}
\label{section:discussion}
\subsection{The case of infinite graphs}
\label{subsec:infinite}
The graph product of groups can be defined for any collection of groups $\G=\{G_i:i\in V\}$ and any simple graph $\Gamma\subset\binom{V}{2}$. However, throughout the paper we assumed that $V=[m]$ is a finite set. We claim that Theorem \ref{thm:explicit_presentation} holds without this assumption. The required changes in the statement and the proof are listed below.
\begin{itemize}
\item We fix a linear order on $V$. In all formulas, $J\subset V$ is assumed to be finite. Hence the elements $\max(J)\in V$ are well defined and $\K_J$ are finite simplicial complexes.
\item In the definition of simplicial complexes (and hence in the definition of the clique complex $\K(\Gamma)$), we require that all faces $I\in\K$ are finite.
\item As in \cite{graph_wreath} and \cite{davis_asphericity_corrigenda}, infinite polyhedral products of pointed topological pairs $(\underline{X},\underline{A})=\{(X_i,A_i,\ast_i),i\in V\}$ are defined as unions
$$(\underline{X},\underline{A})^\K:=\bigcup_{I\in\K}\bigcup_{\begin{smallmatrix}
J\subset V\setminus I,\\|J|<\infty\end{smallmatrix}}
\Big(
\prod_{i\in I}X_i\times
\prod_{j\in J}A_j\times
\prod_{i\notin I\sqcup J}\ast_i\Big)
\subset\prod_{i\in V}X_i,$$ and are endowed with the colimit topology. Equivalently, consider the \emph{restricted direct product} $$\bigoplus_{i\in I}X_i:=\{\underline{x}\in\prod_{i\in I}X_i:~\{i\in I:x_i\neq\ast_i\}\text{ is finite}\}.$$ Then
we define $(\underline{X},\underline{A})^\K:=\bigcup_{I\in\K}\Big(\prod_{i\in I}X_i\times\bigoplus_{j\in V\setminus I}A_j\Big)\subset\bigoplus_{i\in V}X_i.$
\item As noted in \cite{davis_asphericity_corrigenda}, the Davis' approach to the proof of Proposition \ref{prp:polyhedral_products_are_classifying_spaces} (see \cite[Theorem 2.18 + Theorem 2.22]{davis_asphericity}) works in the infinite case. (The categorical tools from \cite{prv}, used by Panov and Veryovkin, also may be applicable here.)
\item In all arguments of Section \ref{section:upper_bound}, replace the direct product $\prod_{i=1}^mG_i$ with the restricted direct product $\bigoplus_{i\in V}G_i=\{\g\in\prod_{i\in V}G_i:\supp\g\text{ is finite}\}$. The paths $p(I)$ are then defined as infinite concatenations, in which only finite number of paths are non-constant. After omitting the constant subpaths, the concatenation is well defined.
\end{itemize}

\subsection{Number of relations in Cartesian subgroups of graph products}
\begin{rmk}
\label{rmk:m_small}

Denote $\pi:=\Bigast_{J\subset[m]}\Pi_1(\K_J)^{\ast n_J},$ so that $\pi_{\ab}=\bigoplus_{J\subset[m]}H_1(\K_J;\ZZ)^{\oplus n_J}.$ The two bounds in Theorems~\ref{thm:presentation_exists_intro} and \ref{thm:lower_bound} can be stated as
$$\rank\pi_{\ab}\leq
\defi\Cart(\G,\K)+\rank\Cart(\G,\K)
\leq\rank\pi;$$
they disagree if and only if 
$\rank \pi_{\ab}<\rank \pi.$ This is possible even if all the groups $\Pi_1(\K_J)$ are abelian: for example, $1=\rank(\ZZ_2\oplus \ZZ_3)<\rank(\ZZ_2\ast\ZZ_3)=2.$

\begin{figure}[h!]
	\centering
	\begin{tikzpicture}[scale=2.5]
\coordinate(1_1) at (  90:1);
\coordinate(2_1) at (  60:1);
\coordinate(3_1) at (  30:1);
\coordinate(4_1) at (   0:1);
\coordinate(1_2) at (- 30:1);
\coordinate(2_2) at (- 60:1);
\coordinate(3_2) at (- 90:1);
\coordinate(4_2) at (-120:1);
\coordinate(1_3) at (-150:1);
\coordinate(2_3) at (-180:1);
\coordinate(3_3) at ( 150:1);
\coordinate(4_3) at ( 120:1);

\coordinate(5)  at (  90:0.5);
\coordinate(6)  at (  50:0.5);
\coordinate(7)  at (  10:0.5);
\coordinate(8)  at (- 30:0.5);
\coordinate(9)  at (- 70:0.5);
\coordinate(10) at (-110:0.5);
\coordinate(11) at (-150:0.5);
\coordinate(12) at ( 170:0.5);
\coordinate(13) at ( 130:0.5);

\draw (1_1)--(2_1)--(3_1)--(4_1)--(1_2)--(2_2)--(3_2)--(4_2)--(1_3)--(2_3)--(3_3)--(4_3)--(1_1);

\draw (5)--(6)--(7)--(8)--(9)--(10)--(11)--(12)--(13)--(5);
\draw (5)--(10)--(6)--(11)--(7)--(12)--(8)--(13)--(9)--(5);
\draw (1_1)--(5);
\draw (1_2)--(8);
\draw (1_3)--(11);
\draw (5)--(2_1)--(6)--(3_1)--(7)--(4_1)--(8)--(2_2)--(9)--(3_2)--(10)--(4_2)--(11)--(2_3)--(12)--(3_3)--(13)--(4_3)--(5);

\fill (1_1) circle (1pt);
\fill (2_1) circle (1pt);
\fill (3_1) circle (1pt);
\fill (4_1) circle (1pt);
\fill (1_2) circle (1pt);
\fill (2_2) circle (1pt);
\fill (3_2) circle (1pt);
\fill (4_2) circle (1pt);
\fill (1_3) circle (1pt);
\fill (2_3) circle (1pt);
\fill (3_3) circle (1pt);
\fill (4_3) circle (1pt);
\fill (5) circle  (1pt);
\fill (6) circle  (1pt);
\fill (7) circle  (1pt);
\fill (8) circle  (1pt);
\fill (9) circle  (1pt);
\fill (10) circle (1pt);
\fill (11) circle (1pt);
\fill (12) circle (1pt);
\fill (13) circle (1pt);

\draw (1_1) node[above] {$1$};
\draw (2_1) node[above right] {$2$};
\draw (3_1) node[above right] {$3$};
\draw (4_1) node[right] {$4$};
\draw (1_2) node[below right] {$1$};
\draw (2_2) node[below right] {$2$};
\draw (3_2) node[below] {$3$};
\draw (4_2) node[below left] {$4$};
\draw (1_3) node[left] {$1$};
\draw (2_3) node[left] {$2$};
\draw (3_3) node[above left] {$3$};
\draw (4_3) node[above] {$4$};

\draw (5) node[above] {$5$};
\draw (6) node[above right] {$6$};
\draw (7) node[above] {$7$};
\draw (8) node[below right] {$8$};
\draw (9) node[below] {$9$};
\draw (10) node[below] {$10$};
\draw (11) node[below left] {$11$};
\draw (12) node[above left] {$12$};
\draw (13) node[above left] {$13$};
		\end{tikzpicture}
	\caption{The 1-skeleton of a flag complex with $\rank\pi_{\ab} < \rank\pi.$
	}
	\label{fig:min_example}
\end{figure}
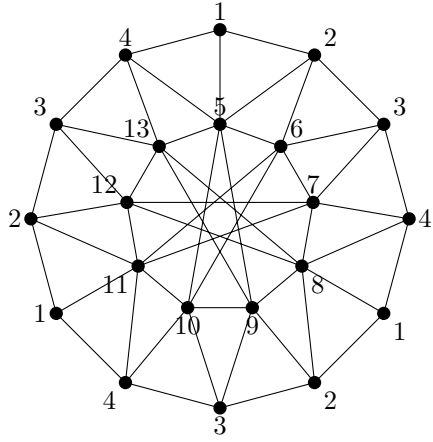
The smallest flag complex with $\rank\pi>\rank\pi_{\ab}$ known to the author is the clique complex of the graph from~Figure~\ref{fig:min_example}. This two-dimensional complex has $13$ vertices and its fundamental group is the \emph{trefoil knot} group $H=\langle a,b\mid a^2=b^3\rangle$. Since $\rank H=2$, $\rank H_{\ab}=1$ and $H=\Pi_1(\K)$ is a free summand of $\pi,$ it follows that $\rank\pi>\rank\pi_{\ab}.$
\end{rmk}
\begin{prb}
Find the deficiency of $\RC'_\K$ for $\K$ as in Figure \ref{fig:min_example}.
\end{prb}

Although $\rank\pi$ and $\rank\pi_{\ab}$ can differ significantly, they agree when $m$ is small enough. Also, these ranks agree for surface triangulations:
\begin{prp}
\label{prp:flag_surface_triangulation}
Let $\K$ be a full subcomplex in a flag triangulation of a closed oriented surface. Then
$\rank\bigoplus_{J\subset [m]} H_1(\K_J;\ZZ)^{\oplus n_J} = \sum_{J\subset[m]}n_J\rank\Pi_1(\K_J).$
\end{prp}
\begin{proof}
It is sufficient to show that, for each connected full subcomplex $\K_J$ of $\K$, the group $H_1(\K_J;\ZZ)$ is free abelian, and $\rank H_1(\K_J;\ZZ)=\rank\pi_1(\K_J)$.

Indeed: for $|\K_J|$ a closed oriented surface, this is well known. Otherwise $|\K_J|$ is a proper compact subset of a closed surface with finitely generated fundamental group, hence $\pi_1(\K_J)$ is free by \cite[Corollary 7]{fischer-zastrow}. (More elementary, one can prove that $\pi_1(\K_J)$ is free by collapsing $\K_J$ onto its 1-skeleton.)
\end{proof}

Due to Proposition \ref{prp:polyhedral_products_are_classifying_spaces} and \cite[Theorem 10.6.1]{davis} (see also \cite[Theorem 2.3]{licai_products}), $\RK$ is an aspherical $3$-manifold if and only if $\K$ is a flag triangulation of $S^2.$ By \cite[Theorem 2.5]{epstein}, fundamental group of any aspherical $3$-manifold admits a \emph{balanced} presentation, i.e. a presentation by $n$ generators modulo $n$ relations for some $n$. For $\RK$, our methods give a \emph{minimal} balanced presentation.
\begin{prp}
Let $\K$ be a flag triangulation of $S^2$. Then Theorem \ref{thm:explicit_presentation} provides a presentation of the $3$-manifold group $\RC'_\K$ by $N(\K)$ generators modulo $N(\K)$ relations, where $N(\K)=\rank\RC'_\K=\sum_{J\subset[m]}\widetilde{b}_0(\K_J)$.
\end{prp}
\begin{proof}
The group $H_1(\RK;\ZZ)\cong\bigoplus_{J\subset[m]}\H_0(\K_J;\ZZ)\simeq\ZZ^{N(\K)}$ is free abelian, so $\bigoplus_{J\subset[m]}\H_1(\K_J;\ZZ)\cong H_2(\RK;\ZZ)\simeq \ZZ^{N(\K)}$ by the Poincar\'e duality. We have $\sum_{J\subset[m]}\rank\Pi_1(\K_J)=N(\K)$ by Proposition \ref{prp:flag_surface_triangulation}. Hence the presentation from Corollary \ref{crl:presentations_for_rck} is minimal and balanced.
\end{proof}

It would be interesting to check if the geometric methods of Gruji\'c \cite{grujic} give a similar presentation.

\subsection{Finite coverings and quotients of real moment-angle complexes}
The group $\{\pm 1\}^m$ acts naturally on the real moment-angle complex $\RK=([-1,1],\{-1,1\})^\K$. Quotients $\RK/H$ by freely acting subgroups $H\subset\{\pm 1\}^m$ are known as \emph{real partial quotients}. In this way one can construct manifolds which admit interesting geometric structures \cite{wu_yu,erokhovets}; the most important examples are the \emph{small covers}, the real analogues of quasitoric manifolds \cite{dj}. The fundamental group of a real partial quotients $\pi_1(\RK/H)$ fits into the short exact sequence
$$1\to \pi_1(\RK/H)\to \RC_\K\to \{\pm 1\}^m/H\to 1;$$
in particular, $\RC'_\K\subset\pi_1(\RK/H)\subset\RC_\K$. Classifying spaces of similar subgroups in the right-angled Artin groups were considered in \cite{gilv}.

Explicit presentations for fundamental groups of small covers were given by Wu and Yu \cite{wu_yu}; in the three-dimensional case, Gruji\'c \cite{grujic} described a minimal balanced presentation.

Note that $\ZZ[1/2]$-homology of real toric spaces \cite{choi_park} and $\ZZ$-homology of small covers \cite{cai_choi} is known. Also, $\RK/H$ decomposes as a union of $2^m/|H|$ copies of $\ccK$ by construction. Hence one can argue as in Sections \ref{section:lower_bound} and \ref{section:upper_bound} to obtain small presentations of these groups as well as bounds on their rank and deficiency.

\begin{prb} Study the fundamental groups of real partial quotients using the methods of this paper.
\end{prb}

The same methods apply for fundamental groups of covering spaces for $\RK$, which are identified with subgroups in $\RC'_\K$.

\subsection{Similarities with graph products of algebras}
\label{subsection:comparison}
Our calculation of relations in the group $\RC_\K'=\pi_1(\RK)$ continues the list of parallel results between the ``real'' and ``complex'' moment-angle complexes $\RK=(D^1,S^0)^\K\simeq (\mathbb{R},\mathbb{R}\setminus\{0\})^\K$ and $\ZK:=(D^2,S^1)^\K\simeq(\mathbb{C},\mathbb{C}\setminus\{0\})^\K.$ Many such results are discussed in \cite[Introduction + Section 7]{prah}. Below we expand this list.

Here $\K$ is a flag simplicial complex on $[m],$ and $\k$ is a principal ideal domain. Note that for flag complexes $\RK$ is aspherical (hence its homotopy properties are governed by the group $\pi_1(\RK)$), while $\ZK$ is simply connected and coformal \cite[Corollary 6.8]{hozk_flag} (hence its homotopy properties are closely related to the associative algebra $H_*(\Omega\ZK;\k)$).
  
On the ``associative algebra side'', we have:
\begin{enumerate}
\item $\Tor^{H_*(\Omega\ZK;\k)}_i(\k,\k)\cong\bigoplus_{J\subset[m]}\H_{i-1}(\K_J;\k)$ by \cite[Theorem 1.2]{cat(zk)};
\item $\mathrm{cat}(\ZK)=1+\max_{J\subset[m]}\mathrm{cdim}_\ZZ\,\K_J$ by \cite[Theorem 1.3]{cat(zk)}, where $\mathrm{cat}(X)$ is the Lusternik--Schnirelmann category of $X$;
\item The $\ZZ\times\ZZ_{\geq 0}^m$-graded $\k$-algebra $H_*(\Omega\ZK;\k)$ can be presented by
$$\sum_{J\subset[m]}\b_0(\K_J)\text{ generators modulo }
\sum_{J\subset[m]}\mathrm{gen}\,H_1(\K_J;\k)\text{  relations}$$
\cite[Theorem 5.5]{hozk_flag}, where the generators and relations are homogeneous; this presentation is minimal among $\ZZ\times\ZZ_{\geq 0}^m$-homogeneous presentations. (If we consider $H_*(\Omega\ZK;\k)$ as a $\ZZ$-graded algebra, we actually need $\sum_{n\geq 0}\mathrm{gen}\,\bigoplus_{|J|=n}H_1(\K_J;\k)$ homogeneous relations \cite[Theorem 5.6]{hozk_flag}).
\item More generally, if $\underline{X}=(X_1,\dots,X_m)$ are simply connected and $\k$ is a field, then the graded $\k$-algebra $H_*(\Omega(C\Omega\underline{X},\Omega\underline{X})^\K;\k)$ can be presented by
\begin{align*}
&\sum_{J\subset[m]}\b_0(\K_J)\cdot\prod_{i\in J}\dim\H_*(\Omega X_i;\k)\text{ homogeneous generators modulo }\\
&\sum_{J\subset[m]}\b_1(\K_J)\cdot\prod_{i\in J}\dim\H_*(\Omega X_i;\k)\text{ homogeneous relations}
\end{align*}
\cite[Corollary 5.14 + Proposition 4.4]{licai_simplicial}.
\end{enumerate}

For a discrete group $G$, we have $G\simeq\Omega BG$ as topological monoids, hence
$H_*(\Omega BG;\ZZ) = H_0(\Omega BG;\ZZ)\cong\ZZ[G]$ as algebras
and $\mathrm{H}_*(G;\ZZ)=\Tor^{\ZZ[G]}_*(\ZZ,\ZZ)\cong\Tor^{H_*(\Omega BG;\ZZ)}_*(\ZZ,\ZZ)$
as graded abelian groups. Hence on the ``group theory side'' we have
\begin{enumerate}
\item $\Tor^{H_*(\Omega\RK;\ZZ)}_i(\ZZ,\ZZ)\cong\bigoplus_{J\subset[m]}\H_{i-1}(\K_J;\ZZ)$ by Proposition \ref{prp:homology_of_cartgk} and Corollary \ref{crl:RCK_RCK'_classifying};
\item $\mathrm{cat}(\RK)=1+\max_{J\subset[m]}\mathrm{cdim}_\ZZ\,\K_J$ by \cite[Proposition 5.12]{cat(zk)};
\item The group $\pi_1(\RK)=\RC_\K'$ can be presented by
$$\sum_{J\subset[m]}\b_0(\K_J)\text{ generators modulo }\sum_{J\subset[m]}\rank\Pi_1(\K_J)\text{ relations}$$ (Corollary \ref{crl:presentations_for_rck}). The number of generators is minimal.
\item More generally, if $\G=(G_1,\dots,G_m)$ are finite groups, then the group $\Cart(\G,\K)=\pi_1((E\G,\G)^\K)=\pi_1((C\Omega B\G,\Omega B\G)^\K)$ can be presented by
\begin{align*}
&\sum_{J\subset[m]}\b_0(\K_J)\cdot\prod_{i\in J}\rank\H_*(\Omega B G_i;\ZZ)\text{ generators modulo }
\\
&\sum_{J\subset[m]}\rank\Pi_1(\K_J)\cdot\prod_{i\in J}\rank\H_*(\Omega B G_i;\ZZ)\text { relations}
\end{align*} by Theorem \ref{thm:presentation_exists_intro}. (Here we use that $|G|-1=\rank\H_*(\Omega BG;\ZZ)$.)
\end{enumerate}
In Remark \ref{rmk:low_est_for_algebras} we also discuss the similarity between lower bounds for presentations of groups and of connected graded algebras.

\begin{prb}
What can be said about the loop homology of the polyhedral product $(C\Omega\underline{X},\Omega\underline{X})^\K$ if $\K$ is a flag simplicial complex and the spaces $X_i$ are neither simply connected nor aspherical?
\end{prb}

\begin{rmk}
The parallel results about $\ZK$ and $\RK$ can be possibly unified in terms of Morel--Voevodsky's motivic homotopy theory, since $\RK=(D^1,S^0)^\K$ and $\ZK=(D^2,S^1)^\K$ are real (resp., complex) Betti realisations of the \emph{motivic moment-angle complex} $\ZK^{\mathbb{A}^1}:=(\mathbb{A}^1,\mathbb{G}_m)^\K,$ a motivic space recently considered by Hornslien \cite{hornslien}. Moreover, the commutator bracket in $\pi_1(\RK)$ and the Whitehead bracket in $\pi_*(\ZK)$ seem to be the Betti realisations of the (non-bilinear) ``motivic Whitehead bracket'' $G_k\times G_\ell\to G_{k+\ell}$ between the non-abelian groups $G_n:=[\Sigma(\mathbb{G}_m)^{\wedge n},\ZK^{\mathbb{A}^1}]$.
\end{rmk}

\end{document}